\documentclass[11pt]{article}
\usepackage{amsmath,amsthm,amssymb,a4wide}

\newcommand{\beq}{\begin{equation}}
\newcommand{\eeq}{\end{equation}}
\newcommand{\beqa}{\begin{eqnarray}}
\newcommand{\eeqa}{\end{eqnarray}}

\newcommand{\al}{\alpha}
\newcommand{\be}{\beta}
\newcommand{\de}{\delta}
\newcommand{\De}{\Delta}
\newcommand{\ga}{\gamma}
\newcommand{\Ga}{\Gamma}

\newcommand{\la}{\lambda}
\newcommand{\La}{\Lambda}
\newcommand{\om}{\omega}
\newcommand{\Om}{\Omega}
\newcommand{\si}{\sigma}

\newcommand{\ze}{\zeta}
\newcommand{\Th}{\Theta}

\newcommand{\ve}{\varepsilon}
\newcommand{\vp}{\varphi}

\newcommand{\cA}{{\mathcal A}}

\newcommand{\cC}{{\mathcal C}}
\newcommand{\cE}{{\mathcal E}}
\newcommand{\cF}{{\mathcal F}}
\newcommand{\cH}{{\mathcal H}}

\newcommand{\cO}{{\mathcal O}}

\newcommand{\OH}{{{\mathcal O}^1\backslash\mathcal H}}
\newcommand{\dH}{{\Gamma_0(d)\backslash\mathcal H}}

\newcommand{\rz}{{\mathbb R}}
\newcommand{\nz}{{\mathbb N}}
\newcommand{\gz}{{\mathbb Z}}
\newcommand{\kz}{{\mathbb C}}

\newcommand{\qz}{{\mathbb Q}}

\DeclareMathOperator{\nm}{N}
\DeclareMathOperator{\tr}{Tr}
\DeclareMathOperator{\im}{Im}
\DeclareMathOperator{\re}{Re}
\DeclareMathOperator{\cl}{cl}
\DeclareMathOperator{\arcosh}{arccosh}

\newcommand{\ndiv}{\!\mbox{$\not|$}}

\newcommand{\rto}{\rightarrow}
\newcommand{\lto}{\longrightarrow}

\numberwithin{equation}{section}

\newtheorem{theorem}{Theorem}[section]
\newtheorem{lemma}[theorem]{Lemma}
\newtheorem{prop}[theorem]{Proposition}
\newtheorem{cor}[theorem]{Corollary}

\theoremstyle{definition}

\theoremstyle{remark}

\begin{document}

\bibliographystyle{abbrv}

\thispagestyle{empty}

\begin{center}
{\LARGE\bf A spectral correspondence for Maa\ss\ waveforms\footnote{
Supported by Deutscher Akademischer 
Austauschdienst and Svenska Institutet}} \\
\vspace*{2cm}
{\large Jens Bolte}\footnote{Electronic address: {\tt 
bol@physik.uni-ulm.de}}\\ \vspace*{5mm}
Abteilung Theoretische Physik\\ 
Universit\"at Ulm, Albert-Einstein-Allee 11\\
D-89069 Ulm\\ Germany \\
\vspace*{1cm}
{\large Stefan Johansson}\footnote{Electronic address: {\tt 
sj@math.chalmers.se}}\\ \vspace*{5mm}
Department of Mathematics\\
Chalmers University of Technology\\ 
and G\"oteborg University\\
S-412 96 G\"oteborg\\
Sweden\\
\end{center}
\vfill

\begin{abstract}
Let $\cO^1$ be a (cocompact) Fuchsian group, given as the group of units 
of norm one in a maximal order $\cO$ in an indefinite quaternion division 
algebra over $\qz$. Using the (classical) Selberg trace formula, we show 
that the eigenvalues of the automorphic Laplacian for $\cO^1$ and their 
multiplicities coincide with the eigenvalues and multiplicities of the 
Laplacian defined on the Maa\ss\ newforms for the Hecke congruence group 
$\Ga_0(d)$, when $d$ is the discriminant of the maximal order $\cO$.
We also show the equality of the traces of certain Hecke operators 
defined on the Laplace eigenspaces for $\cO^1$ and the newforms of
level $d$, respectively.
\end{abstract}

\vspace*{5mm}

\noindent {\small {\bf Key words:} Maa\ss\ waveforms, Selberg trace formula, 
Hecke operators, arithmetic Fuchsian groups \\
\noindent {\bf AMS(MOS) subject classifications:} 11F72, 30F35, 11F12, 11F32}

\newpage

\setcounter{page}{1}

\section{Introduction}
\label{intro}
In several situations correspondences between spaces of automorphic
forms for cocompact and non-cocompact Fuchsian groups are well-known.
From a representation theoretic point of view, the most general such
correspondence is covered by the Jacquet-Langlands correspondence
\cite{JL}. In certain examples it is, however, also possible to
study correspondences by classical means only. Our focus will be
on non-holomorphic automorphic forms (Maa\ss\ waveforms) of weight 
zero and with trivial multiplier, a situation for which classical 
investigations can e.g.\ be found in \cite{Hejhal3}.

To be precise, let $\cO$ be an order in an indefinite rational 
quaternion division algebra so that its group $\cO^1$ of units of
norm one can be considered as a cocompact Fuchsian group. Maa\ss\ 
waveforms for $\cO^1$, i.e.\ eigenfunctions of the automorphic
Laplacian associated with $\cO^1$, can then be lifted to Maa\ss\ 
cusp forms for the Hecke congruence group $\Ga_0(d)$, where $d$
is the (reduced) discriminant of the order $\cO$ \cite{BJ1}.
These lifts are realised as integral transforms with certain (classical)
Siegel theta functions as kernels. As a consequence, theta-lifts preserve
eigenvalues of the hyperbolic Laplacian. In \cite {BJ1},
we also showed that in the case of so-called Eichler orders the eigenvalues
of Hecke operators also remain unchanged. Furthermore in the language of the 
Atkin-Lehner formalism \cite{AL}, theta-lifts of a Hecke basis for 
$L^2(\cO^1\backslash\cH)$ were  found to be newforms of a level 
dividing $d$. Counting theta-lifts and newforms suggested that 
in the case of maximal orders a Hecke basis for $L^2(\cO^1\backslash\cH)$ 
is lifted to a Hecke basis for the newspace of level $d$.

In this paper, we continue to address the question to what extent
theta-lifts are newforms of level $d$. We concentrate on maximal orders
$\cO$ in indefinite rational quaternion division algebras and exploit 
several versions of the (classical) Selberg trace formula \cite{Selberg1}.
Our first result in this direction is summarised in Theorem 
\ref{mainthm1}: The Laplace eigenvalues and their multiplicities for the 
cocompact group $\cO^1$ and for the newforms of level $d$ coincide. This, 
however, does not yet imply that theta-lifts provide isomorphisms between
Laplace (and Hecke) eigenspaces in $L^2(\cO^1\backslash\cH)$ and $L^2
(\Ga_0(d)\backslash\cH)$, respectively, since the theta-lifts are so far 
neither known to be injective, nor to map into the newspace. Our second 
result, as summarised in Theorem \ref{mainthm2}, states that the traces 
of the Hecke operators $\widetilde T_p$ defined on the Laplace eigenspace
in $L^2(\cO^1\backslash\cH)$ coincide with those of the Hecke operators 
$T_p$ defined on the corresponding eigenspaces of newforms in
$L^2(\Ga_0(d)\backslash 
\cH)$ for the infinitely many primes $p$ not dividing the level $d$. Since
theta-lifts preserve Hecke eigenvalues, one thus concludes that all one 
dimensional Laplace eigenspaces in $L^2(\cO^1\backslash\cH)$ lift to the 
newspace of level $d$.

The outline of this paper is as follows: In Section~\ref{sec2}, we fix our
notation and recall some basic facts about the Fuchsian groups $\cO^1$ 
and $\Ga_0(d)$ which are used in the later sections. We then recall the 
results of our preceding paper \cite{BJ1} that are relevant for the 
present work.
Some arithmetic lemmas that are needed later on are established  
in Section~\ref{sec3}. 
In Section~\ref{sec4}, we first recall several versions of the Selberg 
trace formula and rewrite them in a manner that is suitable for our
later purposes. Then we establish the Selberg trace formula for the
groups $\Ga_0(m)$, $m$ square free, that includes the eigenvalues of 
Hecke operators $T_p$, $p\ndiv m$.
Finally our main results are presented in Section \ref{sec5}. 

This work is a result of a cooperation that was enabled by Deutscher
Akademischer Austauschdienst and Svenska Institutet. We want to thank
both institutions for their support. J.B.\ would like to thank Dennis
Hejhal for a discussion on the use of the trace formula in the context
of this paper.

\section{Background and notations}
\label{sec2}

In this section we fix the notation and also recall some basic
facts. Good references for this well-known material are \cite{Iwaniec}
and \cite{Miyake}, see also \cite{BJ1}.

Throughout the paper $\cO$ will denote an order in an indefinite
rational quaternion algebra $\cA$. Except in Section~\ref{sec3}, $\cA$
will be a division algebra and $\cO$ a maximal order in $\cA$. The
discriminant of $\cO$ will be denoted by $d(\cO)$. If $\cO$ is
maximal, then $d(\cO)$ is equal to the discriminant of $\cA$ and
hence equal to the product of the {\em even} number of primes ramified
in $\cA$. We define 
\[ \cO^{n}=\{ x\in\cO : \nm(x)=n\} \]
for $n\in\gz$, where $\nm:\cO\longrightarrow\gz$ is the (reduced)
norm. For the following, we choose a fixed embedding $\si:\cO
\longrightarrow M_2(\rz)$. It is then well-known that $\si(
\cO^1)\subset SL_2 (\rz)$ is a Fuchsian group, which we also denote by 
$\cO^1$. This group acts on the complex upper half-plane $\cH$ by 
M\"obius transformations so that the orbit space $X_{\cO}=\OH$ can be 
considered as a Riemann surface. This surface is compact, since
$\cA$ is a division algebra. We fix a suitable fundamental domain 
$\cF_{\cO}\subset\cH$ for $\cO^1$ with (hyperbolic) area  $A_{\cO}$. If 
$\cO$ is a maximal order, then
\begin{equation}
  \label{area1}
  A_{\cO}=\frac\pi3\prod_{p|d(\cO)}(p-1). 
\end{equation}

Let $L^{2}(X_{\cO})$ be the usual Hilbert space of functions which are
square integrable on $X_{\cO}$ with respect to the hyperbolic volume 
form. Since $X_\cO$ is compact, the spectrum of the hyperbolic Laplacian 
$-\Delta$ on $L^{2}(X_{\cO})$ is discrete, and is comprised of the 
eigenvalues
\[ 0=\lambda_{0}<\lambda_{1}\leq\lambda_{2}\leq\ldots \;,\quad
\lambda_{n}\lto \infty.\]
The counting function for the eigenvalues (with multiplicities) will be 
denoted by 
\[ N_{\cO}(\lambda)=\#\{n :\lambda_{n}\leq\lambda\},\]
 and 
$\{\vp_{k} : k\in\nz_{0}\}$ will be an orthonormal basis of $L^{2}
(X_{\cO})$ with $-\Delta\vp_{k}=\lambda_{k}\vp_{k}$. Thus $\vp_0$ is a
constant function and spans the eigenspace corresponding to $\la_0=0$,
which we denote by $\kz$ for simplicity. We also introduce 
$L^{2}_0(X_{\cO})$ as the orthogonal complement to $\kz$ within $L^{2}
(X_{\cO})$.  

According to a well-known procedure, one can introduce Hecke operators 
$\widetilde{T}_{n}$, $n\in\gz$, on $L^{2}(X_{\cO})$ which are
self-adjoint and bounded operators that commute among themselves and also 
with the hyperbolic Laplacian. It is therefore possible to choose
the orthonormal basis $\{\vp_{k} : k\in\nz_{0}\}$ for $L^{2}(X_{\cO})$
to consist of joint eigenfunctions of $-\De$ and $\widetilde{T}_{n}$, 
$n\in\gz$, i.e. the Laplace eigenfunctions $\vp_{k}$ also satisfy 
$\widetilde{T}_{n}\vp_{k}=\tilde{t}_{k}(n)\vp_{k}$ for constants
$\tilde{t}_{k}(n)$. Such a basis is called a Hecke basis.

The non-compact surfaces we will consider are $X_{d}=\dH$, where
\[ \Gamma_{0}(d)=\left\{
\begin{pmatrix}
  \alpha &\beta \\
  \gamma & \delta
\end{pmatrix}\in SL_{2}(\gz) : \gamma\equiv 0\mod{d} \right\} \]
is the Hecke congruence group of level $d$ with $d$ {\em square free}.
The infinity and all rational numbers are cusps for $\Ga_0(d)$. If
$\om(d)$ is the number of prime divisors of $d$, then the number of $\Ga_0
(d)$-equivalence classes of cusps is $2^{\om(d)}$. As a set of 
representatives for these we choose 
$$F(d)=\left\{\frac{1}{v}:\ v|d,\ v>0\right\}.$$ 
A suitable fundamental domain $\cF_{d}$ for $\Ga_0(d)$ adopted to this 
choice then extends to $\rz$ exactly in the points of $F(d)$. The
hyperbolic area 
$A_{d}$ of $\cF_{d}$ satisfies
\begin{equation}
  \label{area2}
  A_{d}=\frac\pi3\prod_{p|d}(p+1).   
\end{equation}
Associated to the cusps $\frac{1}{v}\in F(d)$ are their stability groups 
\[ \Ga_v=\left\{\ga\in\Ga_0(d):\ \ga\frac{1}{v}=\frac{1}{v}\right\}.\] 
These are conjugate to $\Ga_\infty =\{\pm \bigl(
\begin{smallmatrix}1&j\\0&1\end{smallmatrix}
\bigr):\ j\in\gz\}$ via some $\si_v\in SL_2(\rz)$, $\si_v\infty =\frac{1}{v}$,
$\si_v^{-1}\Ga_v\si_v=\Ga_\infty$. For each $\frac{1}{v}\in F(d)$ one 
also defines a non-holomorphic Eisenstein series $E_{1/v}(z,s)$, 
$z\in\cH$, $s\in\kz$, by
\begin{equation}
\label{Eisendef}
E_{1/v}(z,s)=\sum_{\ga\in\Ga_v\backslash\Ga_0(d)}
\left[\im(\si_v^{-1}\ga z)\right]^s\ ,\quad \re s>1\ ,
\end{equation}
which can be continued to a meromorphic function in $s\in\kz$.

According to the Roelcke-Selberg spectral resolution of the hyperbolic
Laplacian on $L^{2}(X_{d})$, one has an orthogonal decomposition
\[ L^{2}(X_{d})=\cE_{d}\oplus\kz\oplus \cC_{d}, \]
where $\cE_{d}$ is spanned by the Eisenstein series (\ref{Eisendef})
analytically continued to $\re s=\frac{1}{2}$, and $\cC_{d}$ is spanned 
by the cusp forms. The spectrum of $-\Delta$, when restricted to $\kz
\oplus\cC_{d}$, is discrete and the eigenvalues satisfy
\[ 0=\mu_{0}<\mu_{1}\leq\mu_{2}\leq\ldots \;,\quad
\mu_{n}\lto \infty.\]
In analogy to the compact case, we put $N_{d}(\mu)=\#\{n : \mu_{n}\leq
\mu\}$ and choose an orthonormal basis $\{g_{k} : k\in\nz_{0}\}$ of 
$\kz\oplus\cC_{d}$ with $-\Delta g_{k}=\mu_{k}g_{k}$. The Hecke operators 
on $L^{2}(X_{d})$ will be denoted by $T_{n}$, $n\in\nz$, and in this case 
it is possible to choose a Hecke basis $\{g_{k}\}$ such that 
$T_{n}g_{k}=t_{k}(n)g_{k}$ for all $n$ coprime to $d$.

In \cite{BJ1}, we constructed bounded linear operators 
\begin{equation}
\label{Thetalift}
\Theta :\ L^{2}_{0}(X_{\cO})\lto \cC_{d(\cO)} \quad \text{and} \quad
\widetilde{\Theta} :\ \cC_{d(\cO)}\lto L^{2}_{0}(X_{\cO}) 
\end{equation}
preserving both Laplace and Hecke eigenvalues. That is, if
$\{\vp_{k}\}$ and $\{g_{k}\}$ are Hecke bases of $L^{2}_{0}(X_{\cO})$
and $\cC_{d(\cO)}$, respectively, then 
\[ -\Delta(\Theta(\vp_{k}))=\lambda_{k}\Theta(\vp_{k}) \quad \text{and}
\quad
-\Delta(\widetilde{\Theta}(g_{k}))=\mu_{k}\widetilde{\Theta}(g_{k}), \]
and for $(n,d(\cO))=1$
\[ T_{n}(\Theta(\vp_{k}))=\tilde{t}_{k}(n)\Theta(\vp_{k}) \quad \text{and}
\quad
\widetilde{T}_{n}(\widetilde{\Theta}(g_{k}))=
t_{k}(n)\widetilde{\Theta}(g_{k}).  \]

If $g_{k}$ is in a Hecke basis of $\cC_{d}$, then it has a Fourier
expansion at the cusp at infinity that reads
\begin{equation}
  \label{fourier1}
  g_{k}(\tau)=\sum_{n\neq0}c_{k}(n)\,\sqrt v\,
  K_{ir}(2\pi|n|v)\,e^{2\pi inu}\ ,
\end{equation}
with $\tau=u+iv\in\cH$, $\mu_k=r^{2}+\frac14$ being the Laplace 
eigenvalue of $g_{k}$, and $K_{s}(z)$ denoting a modified Bessel 
function. If $(n,d)=1$, then the $n$-th Fourier coefficient is related to 
the $n$-th Hecke eigenvalue via $c_k(n)=c_k(1)\,t_k(n)$. In \cite{BJ1} 
we determined the Fourier expansion of $\Theta(\vp_{k})$, when $\vp_{k}$
is part of a Hecke basis of $L^{2}_{0}(X_{\cO})$, to be
\begin{equation}
  \label{fourier2}
\Th(\vp_k)(\tau)= 4\,\vp_k(z_0)\sum_{n=1}^\infty\tilde t_k(n)\,\sqrt{v}
                     \,K_{ir}(2\pi nv)\,\left[e^{2\pi inu}+\om_k\,e^{-2
                     \pi inu}\right]\ .
\end{equation}
Here $\omega_{k}\in\{\pm1\}$, and $z_{0}$ is a reference point involved 
in the construction of $\Theta$ which can be chosen such that $\vp_k(z_0)
\neq 0$ for all eigenfunctions in a given Hecke basis.

We now recall the necessary background and results from \cite{BJ1} on 
newforms. Let $a,m\in\nz$, $m<d$, be such that $am|d$, and take some 
$h\in\cC_m$. The inclusion $\Ga_0(d)\subset\Ga_0(m)$ implies that $\cC_m
\subset\cC_d$ so that $h\in\cC_d$, but also $h^{(a)}\in\cC_d$ where 
$h^{(a)}(\tau)=h(a\tau)$.
The linear span of all such forms $h^{(a)}\in\cC_d$ that derive from
all possible $a,m$ is called the oldspace $\cC_d^{old}$. Its orthogonal
complement within $\cC_d$ is the newspace $\cC_d^{new}$, so that
$\cC_d =\cC_d^{old}\oplus\cC_d^{new}$. One can 
introduce a Hecke basis of $\cC_d$ such that one part of this basis 
spans the oldspace, and the remaining part spans the newspace. A Hecke 
eigenform in the newspace is then called a newform. If $h$ is a
newform in $\cC_{m}$, then $h^{(a)}$ is called an oldform in $\cC_{d}$.

It is trivial to check that all $h^{(a)}$ corresponding to a fixed
$h\in\cC_m$ have the same Laplace eigenvalue. If $h\in\cC_m^{new}$,
then there are $\tau(\frac d m)$ forms $h^{(a)}$ in $\cC_d$ corresponding
to $h$, where $\tau(n)$ is the number of positive divisors of $n$.
Let $\de(m,\la)$ be the dimension of the subspace of $\cC_m$ with Laplace
eigenvalue $\la>0$, and let $\de'(m,\la)$ be the dimension of the 
corresponding subspace of $\cC_m^{new}$. Since 
$\cC_d =\cC_d^{old}\oplus\cC_d^{new}$, these satisfy
\[ \de(d,\la)=\sum_{m|d}\tau\!\left(\tfrac d m\right) \de'(m,\la). \]
Inverting this formula, one gets \cite[(6.7)]{AL}
\begin{equation}
\label{newformsum}
\de'(d,\la)=\sum_{m|d}\be\!\left(\tfrac d m\right) \de(m,\la),
\end{equation}
with 
\begin{equation}
  \label{beta}
  \be(n)=\sum_{k|n}\mu(k)\mu(\tfrac{n}{k}),
\end{equation}
 where $\mu(n)$ is 
the M\"obius function. In particular if $n$ is a product of $r$ distinct 
primes, then $\be(n)=(-2)^{r}$.

If $N_d^{'}(\la)=\#\left\{\mu_k\leq\la:\ g_k\in\cC_d^{new}\right\}$, then
using (\ref{newformsum}) we observed in \cite{BJ1} that 
\begin{equation}
\label{equalasym}
\lim_{\la\rto\infty}\frac{N_{d}^{'}(\la)}{N_{\cO}(\la)}\geq1,
\end{equation}
when $\cO$ is a maximal order with discriminant $d$; equality emerges 
under the hypothesis that the Laplace spectrum of $L^{2}(X_{\cO})$ is
simple. In this paper, we will show the much stronger result that 
there is a one-to-one correspondence including multiplicities between
the Laplace spectra of $L_0^{2}(X_{\cO})$ and $\cC_{d}^{new}$. More
generally, we will prove that the traces of the Hecke operators
$\widetilde{T}_{p}$ and $T_{p}$ on the $\lambda$-eigenspaces of 
$L_0^{2}(X_{\cO})$ and $\cC_{d}^{new}$ coincide for $p\ndiv d$.

\section{Auxiliary results}
\label{sec3}

Let $\cO$ be a so-called Eichler order of square free level $m$ and
discriminant $d\cdot m$, where $d$ is a product of an even number of 
distinct primes. We remark that $(d,m)=1$ is implicit in the definition. 
If $d=1$, then
\[ \cO\cong M(m)=\left\{ 
        \begin{pmatrix}
        a & b \\
        c & d 
        \end{pmatrix}\in M_{2}(\gz) :c\equiv 0 \mod{m} \right\}, \]
and if $m=1$, then $\cO$ is a maximal order in an indefinite quaternion
algebra $\cA$ over $\qz$. Let $\cO^1$ denote the group of elements in 
$\cO$ with (reduced) norm equal to $1$. In particular, $\cO^1\cong 
\Gamma_0(m)$ if $d=1$, and otherwise $\cO^1$ is a cocompact arithmetic 
Fuchsian group.

Let $K$ be a quadratic field extension of $\qz\,$ and $B$ an order in 
$K$. Assume that there is an embedding $\iota : K\lto\cA$. The 
order $B$ is said to be optimally embedded into 
$\cO$ with respect to $\iota$  if $\iota(B)=\cO\cap \iota(K)$. 
Two optimal embeddings, 
$\iota_1$ and $\iota_2$, are conjugate by $\ga$ if 
$\iota_1(B)=\ga\cdot\iota_2(B)\cdot\ga^{-1}$.
The number of optimal embeddings of $B$ modulo conjugation by elements 
in $\cO^1$
is given by the formula 

\begin{equation}
\label{optorders}
 E(B,\cO^1)=\cl(B)\cdot \prod_{p|d}\left(1-\Bigl(\frac B p\Bigr)\right)
\prod_{p|m}\left(1+\Bigl(\frac B p\Bigr)\right),
\end{equation}
where $\cl(B)$ is the class number of $B$ and
\[ \Bigl(\frac B p\Bigr)=\begin{cases}
        1& \text{if $p$ is split in $K$ or $B$ is not a maximal 
                order in $K$,}\\
        -1& \text{if $p$ is unramified in $K$ and $B$ is maximal,} \\
        0& \text{if $p$ is ramified in $K$ and $B$ is maximal.}
                          \end{cases} \]

For imaginary extensions $K$ this is a special
case of \cite[2.5.]{Schneider}, noting that the local embedding
number for a prime $p$ dividing $m$ is $\bigl(1+(\frac B p)\bigr)$
\cite[Th.II.3.2]{Vigneras}. However, in this simple case with a rational ground
field, it is easy to see that the factor multiplying the product of
the local embedding numbers is $\cl(B)$ also for real quadratic
extensions. We remark that the fact that all Eichler orders contain
elements with norm $-1$ \cite[(5.6)]{Stefan} is important when deriving
that the factor is exactly $\cl(B)$.

If $\ga\in\cO$, then the conjugacy class of $\ga$ with respect to 
$\cO^1$ is $\{\ga\}_{\cO^1}=\{\al\ga\al^{-1} : \al\in\cO^{1}\}$. Observe 
that $\tr(\ga)=\tr(\al\ga\al^{-1})$, where $\tr :\cO\longrightarrow\gz$ 
is the (reduced) trace, and $\nm(\ga)=\nm(\al\ga\al^{-1})$. For $t,n\in
\gz$, we now define $E(t,n,\cO^1)$ to be the number of conjugacy classes 
$\{\ga\}_{\cO^1}$ in $\cO$ with $\tr(\ga)=t$ and $\nm(\ga)=n$. If 
$\gamma$ is an arbitrary element with $\tr(\ga)=t$ and $\nm(\ga)=n$, 
then \cite[p.96]{Vigneras}
\begin{equation}
\label{optelements}
E(t,n,\cO^1)=\sum_{B\supseteq\gz[\gamma]} E(B,\cO^1), 
\end{equation}
where the sum is over all orders in $\qz(\gamma)$ containing $\gamma$.

We will only need to consider $E(t,n,\cO^1)$ for $n=1$ or $n=p$ a prime, 
and we may restrict to $t>0$ since $E(t,n,\cO^1)=E(-t,n,\cO^1)$. For 
$n=p$, the case $t=p+1$ is exceptional, since this is the only 
situation where $\qz(\ga)$ is not a field. This implies in particular 
that $E(p+1,p,\cO^{1})=0$ when $\cO$ is an order in a division algebra.

\begin{lemma} 
\label{exclemma}
If $p$ is a prime and $m$ is a product of $\om(m)$ distinct primes, then
\[ E\bigl(p+1,p,\Gamma_0(m)\bigr)=2^{\om(m)}(p-1). \]
Moreover, if  $n_0(v)$ are integers satisfying $n_0(v)
\equiv (p-1)v^{-1}\pmod{\frac m v}$ for positive divisors $v$ of $m$, then
\[ 
\bigcup_{\substack{ v|m \\ v>0}}\left\{\ga_{v,n}=\begin{pmatrix}
        p-nv& n\\
        v(p-nv-1)& nv+1   \end{pmatrix} :\ n=n_0(v)+k\frac{m}{v},\;
        0\le k<p-1 \right\} 
\]
is a set of representatives of the different conjugacy classes with trace
$p+1$ and norm $p$.
\end{lemma}

\begin{proof}
Let $\al=\bigl(\begin{smallmatrix}
        a& b\\
        c& d  \end{smallmatrix}\bigr) \in M(m)$ be an arbitrary element 
with $\tr(\al)=p+1$ and $\nm(\al)=p$. Then $\al$ has two rational fixed 
points, namely $\frac{a-1}c$ and $\frac{1-d}c$. Each of these fixed 
points is mapped to a unique member in $F(m)=\{\frac1v : v|m\}$ by
        some elements $\be_1$  
and $\be_2$ in $\Ga_0(m)$, since all rational numbers are cusps for 
$\Ga_0(m)$. Hence we may assume that $\al$ has a fixed point in $F(m)$, 
since we may conjugate by $\be_1$ or $\be_2$. 

Fix a positive divisor $v$ of $m$. 
The element $\ga=\bigl(\begin{smallmatrix}
        1& 0\\
        -v& 1  \end{smallmatrix}\bigr)\in\Ga_0(1)$ maps $\frac1v$ to
$\infty$. Elements with trace $p+1$, norm $p$ and one fixed point in 
$\infty$ are either of the form
$\bigl(\begin{smallmatrix}
        1& n\\
        0& p  \end{smallmatrix}\bigr)$ or
$\bigl(\begin{smallmatrix}
        p& n\\
        0& 1  \end{smallmatrix}\bigr)$.
Hence $\al$ is of one of the forms:
\begin{equation} \label{form1}
 \al=\ga^{-1}\begin{pmatrix}
        1& n\\
        0& p   \end{pmatrix} \ga=\begin{pmatrix}
        1-nv& n\\
        v(1-nv-p)& nv+p   \end{pmatrix},
\end{equation}
or
\begin{equation} \label{form2}
\al=\ga^{-1}\begin{pmatrix}
        p& n\\
        0& 1   \end{pmatrix} \ga=\begin{pmatrix}
        p-nv& n\\
        v(p-nv-1)& nv+1   \end{pmatrix}.
\end{equation}
A priori, $n$ need not be an integer. However, the condition $\al\in 
M(m)$ forces $n$ to be an integer, and moreover is equivalent to 
$\frac{m}{v}|(1-nv-p)$ and $\frac{m}{v}|(p-nv-1)$, respectively. Since 
$(\frac{m}{v},v)=1$, this implies that $\al\in M(m)$ is equivalent to $n$
being congruent to a unique class modulo $\frac{m}{v}$.

There are only two possibilities for two elements of the form 
(\ref{form1}) or (\ref{form2}) to be conjugate by elements in $\Ga_0(m)$. 
Either directly by an element in $\Ga_v$ for some $v|m$, or otherwise 
they correspond to the two different choices of fixed point of $\al$ 
combined with an element in $\Ga_v$. First we will consider $\Ga_v$.

One readily checks that $\Ga_v= \ga^{-1}\left< \be \right>\ga$, where 
 $\be=\bigl(\begin{smallmatrix}
        1& \frac m v\\
        0& 1  \end{smallmatrix}\bigr)$. 
A direct computation gives
\[ (\ga^{-1}\be^k\ga)\begin{pmatrix}
        1-nv& n\\
        v(1-nv-p)& nv+p   \end{pmatrix}(\ga^{-1}\be^k\ga)^{-1}=
        \begin{pmatrix}
        1-n'v& n'\\
        v(1-n'v-p)& n'v+p   \end{pmatrix}, \]
where $n'=n+k(p-1)\frac m v$. Since $n$ is congruent to a unique class 
modulo $\frac{m}{v}$, we conclude that we get $p-1$ different conjugacy 
classes. The form (\ref{form2}) is completely analogous. This gives us 
at most $2(p-1)$ different conjugacy classes associated to each of the 
different cusps.

A direct computation shows that conjugation with an element taking the 
fixed point $\frac{a-1}c$ to the cusp $\frac1v\in F(m)$ always gives an 
element of the form (\ref{form2}), and for $\frac{1-d}c$ of the form
(\ref{form1}). Hence we can conclude that the number of different 
conjugacy classes is $2^{\om(m)}(p-1)$. Furthermore, they can be chosen 
to be either of the form (\ref{form1}) or (\ref{form2}).
\end{proof}

The basis for the definition of the Hecke operator $T_p$ with respect to
$\Ga_0(m)$ for $(p,m)=1$ is the set $M_p(m)=\{\al\in M(m) : \det(\al)
=p\}$. It has the following decomposition \cite[(4.5.25)]{Miyake}
\begin{equation}
\label{decomp}
M_p(m)=\bigcup_{j=0}^{p-1}\Ga_0(m)
\begin{pmatrix}
        1&j\\
        0&p
\end{pmatrix}
\cup\Ga_0(m)    
\begin{pmatrix}
        p&0\\
        0&1
\end{pmatrix}.                  
\end{equation}

\begin{lemma}
\label{ineqcusps}
Assume that $p$ is a prime and $m$ is a square free integer with $(p,m)=1$. 
Then two cusps which are inequivalent modulo $\Ga_0(m)$ are also 
inequivalent modulo $M_p(m)$.
\end{lemma}

\begin{proof}
We recall that $\{\frac1v : v|m\}$ is a set of representatives of the 
different cusps modulo $\Ga_0(m)$. From (\ref{decomp}) it follows that it 
suffices to show that $\ga\frac1v$ is $\Ga_0(m)$-equivalent to 
$\frac1v$ for any element $\ga$ of the form 
$\bigl(\begin{smallmatrix}
        1& j\\
        0& p  
\end{smallmatrix}\bigr) $ or
$\bigl(\begin{smallmatrix}
        p& 0\\
        0& 1  
\end{smallmatrix}\bigr) $, that is, $\frac{1+jv}{pv}$ and $\frac p v$ are
$\Ga_0(m)$-equivalent to $\frac1v$. More generally, we show that 
$\frac1v$ is equivalent to $\frac{x}{yv}$ whenever $(x,yv)=1$ and 
$(y,m)=1$. Take 
$\al=\bigl(\begin{smallmatrix}
        a& b\\
        mc& d  
\end{smallmatrix}\bigr)\in\Ga_0(m)$. Then $\al\frac1v=\frac{x}{yv}$ 
and $\al\in\Ga_0(m)$ is equivalent to
\begin{align*}
1&= ad-mbc \\
x&=a+bv \\
y&=d+c\tfrac m v.
\end{align*}
Substituting the expressions of $a$ and $d$ into the first equality, we 
get
\[ 1=(y-\tfrac m v c)(x-bv)-mbc=xy-c\tfrac m v x-bvy. \]
This has a solution $b,c$, since $(\tfrac m v x,vy)=1$.
\end{proof}

The next three lemmas contain identities which are the main
ingredients when comparing the trace formulae in Section~\ref{sec5}.

\begin{lemma}
\label{embeddinglemma}
Let $\cO$ be a maximal order in an indefinite rational quaternion algebra
with discriminant $d$. Then
\begin{equation}
\label{sumid1}
E(B,\cO^1)=\sum_{m|d}\be\!\left(\tfrac d m\right) E\bigl(B,\Ga_0(m) \bigr)
\end{equation}
and
\begin{equation}
\label{sumid2}
E(t,n,\cO^1)=\sum_{m|d}\be\!\left(\tfrac d m\right)
E\bigl(t,n,\Ga_0(m) \bigr), 
\end{equation}
where $\be$ is given by $(\ref{beta})$.
\end{lemma}

\begin{proof}
The identity (\ref{sumid2}) follows directly from (\ref{sumid1}) by
using (\ref{optelements}). Since $\cl(B)$ only depends on $B$, it follows
from (\ref{optorders}) that (\ref{sumid1}) is equivalent to
\begin{equation}
\label{embeddingid}
\prod_{p|d}\left(1-\Bigl(\frac B p\Bigr)\right)=\sum_{m|d}
\be\!\left(\tfrac d m\right)\prod_{p|m}\left(1+\Bigl(\frac B p\Bigr)
\right).
\end{equation}

That $\cO$ is a maximal order in an indefinite quaternion algebra implies
that $d$ is a product of an even number, $2r$, of different primes. Fix
an order $B$ and let $k=\#\{p : p|d,\;(\frac B p)=1\}$ and
$e=\#\{p : p|d,\;(\frac B p)=0\}$.
We obviously have
\[
\prod_{p|d}\left(1-\Bigl(\frac B p\Bigr)\right)=
        \begin{cases}
        0& \text{if }k>0, \\
        2^{2r-e}& \text{otherwise.}
        \end{cases}
\]

Moreover, a term of the right-hand side of (\ref{embeddingid}) is 
non-zero only for all possible products of the $k+e$ primes with 
$(\frac B p)\neq -1$. This observation yields
\begin{align*}
\sum_{m|d}\be\!\left(\tfrac d m\right)\prod_{p|m}
\left(1+\Bigl(\frac B p\Bigr)\right)&= \sum_{j=0}^e\sum_{i=0}^k 
(-2)^{2r-i-j}\cdot 2^i\cdot 1^j \binom k i\binom e j \\
&= 2^{2r-e}\sum_{j=0}^e (-1)^j 2^{e-j}\binom e j\sum_{i=0}^k (-1)^i\binom k i.
\end{align*}
The desired result follows by observing that 
\[
\sum_{j=0}^e (-1)^j 2^{e-j}\binom e j=(2-1)^e=1 \quad\text{ and }\quad
\sum_{i=0}^k (-1)^i\binom k i=
        \begin{cases}
        (1-1)^k=0& \text{if } k>0, \\
        1& \text{if }k=0.
        \end{cases}
\]
\vspace{-5ex}

\end{proof}

A hyperbolic element $\ga\in\Ga$, $\Ga$ a Fuchsian group, is called 
primitive in $\Ga$, if there is no element $\al\in\Ga$ such that 
$\ga=\al^r$, with $r\ge 2$. For $t>2$, we define $E'(t,1,\cO^1)$ to be 
the number of conjugacy classes in $\cO^1$ of primitive elements with 
trace $t$. We remark that if $\ga$ is hyperbolic, then for any $r\ge1$ 
$\{\ga\}_{\cO^1}=\{\be\}_{\cO^1}$
iff $\{\ga^r\}_{\cO^1}=\{\be^r\}_{\cO^1}$. Moreover, $r\longmapsto|\tr(
\ga^r)|$ is a strictly increasing function for $r\ge1$ when $\ga$ is 
hyperbolic. Using these two facts it is clear that
\begin{equation}
\label{indstep}
E(t,1,\cO^1)=\sideset{}{'}\sum_{s\le t}E'(s,1,\cO^1),
\end{equation}
where the sum is over all $s$ such that if $\tr(\ga)=s$, then there is an 
$r\ge1$ such that $\tr(\ga^r)=t$.
Since $E(3,1,\cO^1)=E'(3,1,\cO^1)$, the following lemma follows by 
straightforward induction on $t$ using (\ref{sumid2}) and (\ref{indstep}).

\begin{lemma}
\label{primitive}
Let $\cO$ be a maximal order in an indefinite rational quaternion algebra
with discriminant $d$. Then
\[
E'(t,1,\cO^1)=\sum_{m|d}\be\!\left(\tfrac d m\right) E'\bigl(t,1,\Ga_0(m) 
\bigr).
\]
\end{lemma}

We conclude this section with three simple arithmetic identities which are
important in the proof of the main theorems.

\begin{lemma}
\label{arithid}
If $d>1$ is a  square free integer and $\om(m)$ the number of primes
dividing $m$, then
\[ \sum_{m|d}\be\!\left(\tfrac d m\right)=(-1)^{\omega(d)} \quad\text{
  and }\quad  
\sum_{m|d}\be\!\left(\tfrac d m\right)2^{\om(m)}=0. \]
Moreover, if $f$ is an arbitrary function and $d$ has at least two
prime divisors, then
\[ \sum_{m|d}\be\!\left(\tfrac d m\right)2^{\om(m)}\sum_{p|m}f(p)=0, \]
where the first sum is over all divisors and the second only over prime 
divisors.
\end{lemma}

\begin{proof}
The first identity follows from the binomial theorem, since
\[ \sum_{m|d}\be\!\left(\tfrac d m\right)=\sum_{i=0}^{\om(d)}(-2)^{\om(d)
-i}\,\binom{\om(d)}i=(1-2)^{\omega(d)}=(-1)^{\omega(d)}. \]
The second identity follows analogously,
\[ \sum_{m|d}\be\!\left(\tfrac d m\right)2^{\om(m)}= 
\sum_{i=0}^{\om(d)}(-2)^{\om(d)-i}\cdot 2^i\binom{\om(d)}i=
(2-2)^{\om(d)}=0. \]
>From this identity we derive the third one by
\begin{align*}
\sum_{m|d}\be\!\left(\tfrac d m\right)2^{\om(m)}\sum_{p|m}f(p)&= 
\sum_{p|d}f(p)\sum_{m|\frac d p}\be\!\left(\tfrac d{mp}\right)2^{\om(mp)} \\
&= 2\cdot\sum_{p|d}f(p)\sum_{m|\frac d p}\be\!\left(\tfrac{d/p}{m}\right)
2^{\om(m)}=0.
\end{align*}
\vspace{-4.5ex}

\end{proof}

\section{Trace formulae}
\label{sec4}
The basic tool we are going to apply in order to prove our main result
will be the Selberg trace formula. We will need several versions of it;
with and without inclusion of Hecke-eigenvalues, for the cocompact
groups $\cO^1$ and the Hecke congruence groups $\Ga_0(m)$. 
Without inclusion of Hecke-eigenvalues the trace formula for all 
Fuchsian groups occurring are well-known. When considering Hecke 
operators $\widetilde T_p$, $T_p$, for primes $p\ndiv d(\cO)$, the 
relevant trace formula for the groups $\Ga_0(m)$, $m>1$, still
has to be evaluated explicitly. The result for $\cO^1$ can be found
in \cite[ch.V,Thm.8.1]{Hejhal1} and for $\Ga_0(1)$ in \cite[(11.10)]
{HejhalDuke}. In this section, we will first recall the known results,
and thereby take the opportunity to introduce some notation that will 
be useful. Then we will derive the necessary trace 
formula not yet to be found in the literature.

In the sequel $h:\,\kz\rto\kz$ always denotes a function  
satisfying
\begin{itemize}
\item $h(r)=h(-r)$,
\item $h(r)$ is holomorphic in the strip $|\im r|\leq\frac{1}{2}+\ve$,
for some $\ve >0$,
\item $|h(r)|\leq C\,(1+\re r)^{-2-\de}$ for some $C>0$ and $\de>0$.
\end{itemize}
The Fourier transform of $h$ will then be written as
\begin{equation*}
\hat h(u)=\frac{1}{2\pi}\int_{-\infty}^{+\infty}h(r)\,e^{-iru}\ dr\ .
\end{equation*}
Since the unit group $\cO^1$ is a cocompact Fuchsian group, the Selberg 
trace formula without Hecke-eigenvalues reads as follows:

\begin{prop}
\label{STFcomp}
Let $\la_k =r_k^2+\frac{1}{4}$ run through all eigenvalues of the
hyperbolic Laplacian on $L^2(X_\cO)$, counted with multiplicities.
Then
\begin{equation*}
\begin{split}
\sum_{k=0}^\infty h(r_k)= 
     &\frac{A_\cO}{4\pi}\int_{-\infty}^{+\infty}h(r)\,r\,\tanh
      (\pi r)\ dr \\
     &+\sum_{t\in\{0,1\}}\frac{E'(t,1,\cO^1)}{2m_t}\sum_{k=1}^{m_t -1}
      \frac{1}{\sin(\frac{k\pi}{m_t})}\int_{-\infty}^{+\infty}h(r)\,
      \frac{e^{-\frac{2k\pi r}{m_t}}}{1+e^{-2\pi r}}\ dr \\
     &+\sum_{t=3}^\infty E'(t,1,\cO^1)\,\arcosh(\tfrac{t}{2})
      \sum_{k=1}^\infty\frac{\hat h\left(2k\arcosh(\frac{t}{2})
      \right)}{\sinh\left(k\arcosh(\frac{t}{2})\right)}\ .\\
\end{split}
\end{equation*}
\end{prop}

\begin{proof}
Recall the trace formula for cocompact Fuchsian groups \cite[(3.2)]
{Selberg1}, \cite[ch.II,Thm.5.1]{Hejhal1}. The sums over representatives
of the primitive elliptic and hyperbolic conjugacy classes of
elements in $\cO^1$ are rewritten as sums over the
traces $t\in\nz_0$ of the representatives. The number of primitive 
conjugacy classes with trace $t$ is $E'(t,1,\cO^1)$. In the elliptic
case, where $t\in\{0,1\}$, $m_t$ denotes 
the order of the primitive element with trace $t$. 
\end{proof}

We now consider the Hecke operators $\widetilde T_p$ when $p$ is
a prime not dividing $d(\cO)$. In this case the relevant trace 
formula is given by

\begin{prop}
\label{STFTpO}
Let $\la_k =r_k^2+\frac{1}{4}$ run through all eigenvalues of the
hyperbolic Laplacian on $L^2(X_\cO)$, counted with multiplicities.
Fix a prime $p\ndiv d(\cO)$ and denote by $\tilde t_k(p)$ the 
eigenvalues of $\widetilde T_p$ on a Hecke basis $\{\vp_k:\ k\in
\nz_0\}$ of $L^2(X_\cO)$. Then
\begin{equation*}
\begin{split}
\sum_{k=0}^\infty\tilde t_k(p)\,h(r_k)= 
     &\sum_{\substack{\{\ga\}_{\cO^1} \\ \ga\,\text{ellip.}}}
      \frac{1}{m_\ga\sqrt{4p-t^2}}\int_{-\infty}^{+\infty}h(r)\,
      \frac{e^{-2r\arcsin\sqrt{1-\frac{t^2}{4p}}}}{1+e^{-2\pi r}}
      \ dr \\
     &+\frac{1}{\sqrt{p}}\sum_{\substack{\{\ga\}_{\cO^1} \\
     \ga\,\text{hyperb.}}} 
      \frac{\arcosh(|\log\ve_\ga|)\,\hat h\left(2\arcosh(\tfrac{t}
      {2\sqrt{p}})\right)}{2\sinh\left(\arcosh(\tfrac{t}
      {2\sqrt{p}})\right)}\ ,\\
\end{split}
\end{equation*}
where $t=\tr(\ga)$ and $p=\nm(\ga)$. Moreover, $m_\ga$ is the order of the 
centraliser $Z_{\cO^1}(\ga)$ when $\ga\in\cO^p$ is elliptic, and  
$\ve_\ga$ is a generator of $Z_{\cO^1}(\ga)$ when $\ga\in\cO^p$ is
hyperbolic. 
\end{prop}
 
\begin{proof}
In order to define the Hecke operator $\widetilde T_p$,
\begin{equation*}
\widetilde T_p\vp(z)=\frac{1}{\sqrt{p}}\sum_{j=1}^{d(p)}
\vp(\ga_j z)\ ,\ \ \ \ \vp\in L^2(X_\cO)\ ,
\end{equation*}
one needs the decomposition
\begin{equation*}
\cO^p=\left\{ u\in\cO:\ \nm(u)=p\right\}=\bigcup_{j=1}^{d(p)}
\cO^1\ga_j\ ,
\end{equation*}
which derives from the decomposition
\begin{equation*}
\cO^1 u\cO^1 =\bigcup_{j=1}^{d(u)}\cO^1\ga_j
\end{equation*}
of the distinct double cosets $\cO^1 u\cO^1$ with $u\in\cO^p$.

Clearly, $\widetilde T_p$ commutes with $-\De$, and also with any
integral operator 
\begin{equation*}
L\vp(z)=\int_{\cF_\cO} K(z,w)\,\vp(w)\ d\mu(w)\ ,\ \ \ \ \vp\in 
L^2(X_\cO)\ ,
\end{equation*}
which has a kernel of the form
\begin{equation*}
K(z,w)=\sum_{\ga\in\cO^1/\{\pm E_2\}}k(\ga z,w)\ ,\ \ \ \ k(z,w)=
\phi\left(\frac{|z-w|^2}{\im z\im w}\right)\ ,
\end{equation*}
where $d\mu(z)=\frac{dx\,dy}{y^2}$ is the hyperbolic volume form and
$\phi\in C^2_0(\rz)$ is arbitrary. The fact that $\widetilde 
T_p L=L\widetilde T_p$ follows from \cite[ch.V,Prop.2.19]{Hejhal1}
and from the observation that $\widetilde T_p$ is $\frac{1}{\sqrt{p}}$
times a finite sum of operators as in \cite[ch.V,Def.2.9]{Hejhal1}.
According to \cite[ch.V,Prop.2.22]{Hejhal1}, $\widetilde T_p L$ has
an integral kernel
\begin{equation}
\label{Kpsum}
K_p(z,w)=\frac{1}{\sqrt{p}}\sum_{j=1}^{d(p)}\sum_{\ga\in\cO^1/\{\pm 
E_2\}}k(\ga\ga_j z,w)=\frac{1}{\sqrt{p}}\sum_{\ga\in\cO^p/\{\pm E_2\}}
k(\ga z,w)\ .
\end{equation}
Following \cite[ch.V]{Hejhal1} further, one calculates $\tr L\widetilde 
T_p$ on the one hand from (\ref{Kpsum}), and on the other hand from the 
spectral expansion
\begin{equation*}
K_p(z,w)=\sum_{k=0}^\infty\tilde t_k(p)\,h(r_k)\,\vp_k(z)\,
\overline{\vp_k}(w)\ .
\end{equation*} 
This yields
\begin{align}
\label{Tptrgeom}
\sum_{k=0}^\infty\tilde t_k(p)\,h(r_k)
   &=\frac{1}{\sqrt{p}}\sum_{\ga\in\cO^p/\{\pm E_2\}}\int_{\cF_\cO} 
    k(\ga z,z)\ d\mu(z) \nonumber\\
   &=\frac{1}{\sqrt{p}}\sum_{\{\ga\}_{\cO^1}}\sum_{\si\in\cO^1/
    Z_{\cO^1}(\ga)}\int_{\si\cF_\cO}k(\si z,z)\ d\mu(z)\ .
\end{align}
In the second line the first sum extends over all conjugacy classes 
$\{\ga\}_{\cO^1}$ of elements $\ga\in\cO^p/\{\pm E_2\}$. For the 
second sum one needs to know the centraliser
\begin{equation*} 
Z_{\cO^1}(\ga)=\left\{\si\in\cO^1:\ \si\ga=\ga\si\right\}=\qz(\ga)
\cap\cO^1
\end{equation*}
of $\ga\in\cO^p$. Since $\cO$ is a maximal order in a division
algebra, any quadratic  
extension $\qz(\rho)$, $\rho\in\cO$, is a field. Moreover,
$\qz(\rho)\cap\cO$ is the maximal order in $\qz(\rho)$. It 
follows that the centraliser $\qz(\ga)\cap\cO^1$ is the group of units 
of norm one in the maximal order of $\qz(\sqrt{D})$, where $D$ is square 
free such that $t^2-4p =n^2 D$.  Elliptic elements 
of $\cO^p$ are characterised by $t^2 -4p<0$, which implies that $D$ 
is negative. In this case $Z_{\cO^1}(\ga)$ is a finite 
cyclic group of order $m_\ga$. In the hyperbolic case, where 
$t^2 -4p>0$ and hence $D$ is positive, the centraliser is infinite 
cyclic, generated by $\ve_\ga$. Going on as in \cite[ch.V]{Hejhal1}, 
we arrive at the desired result; compare \cite[ch.V,Thm.8.1]{Hejhal1}.
\end{proof}

Next we consider the Hecke congruence groups $\Ga_0(m)$ with $m|d
(\cO)$. Thus $m$ is square free and consist of $\om(m)$ prime divisors.
For the rest of this section the discriminant $d(\cO)$ is not important
so that $m$ denotes an arbitrary square free positive integer. Without 
the inclusion of Hecke-eigenvalues the trace formula for $\Ga_0(m)$
is well-known. We recall \cite[Thm.9.9]{HejhalDuke} together with
\cite[(10.2),(10.4)]{HejhalDuke} and use the same notation as in
Proposition \ref{STFcomp}:

\begin{prop}
\label{STFGam}
Let $\mu_k =r_k^2+\frac{1}{4}$ run through all eigenvalues of the
hyperbolic Laplacian on $L^2(X_m)$, counted with multiplicities.
Then
\begin{equation*}
\begin{split}
\sum_{k=0}^\infty h(r_k)= 
     &\frac{A_m}{4\pi}\int_{-\infty}^{+\infty}h(r)\,r\,\tanh
      (\pi r)\ dr\\
     &+\sum_{t\in\{0,1\}}\frac{E'(t,1,\Ga_0(m))}{2m_t}\sum_{k=1}^{m_t 
      -1}\frac{1}{\sin(\frac{k\pi}{m_t})}\int_{-\infty}^{+\infty}h(r)
      \,\frac{e^{-\frac{2k\pi r}{m_t}}}{1+e^{-2\pi r}}\ dr\\
     &+\sum_{t=3}^\infty E'(t,1,\Ga_0(m))\,\arcosh(\tfrac{t}{2})
      \sum_{k=1}^\infty\frac{\hat h\left(2k\arcosh(\frac{t}{2})
      \right)}{\sinh\left(k\arcosh(\frac{t}{2})\right)}\\
     &+2^{\om(m)}\biggl\{\hat h(0)\log(\tfrac{\pi}{2})-\frac{1}{2\pi}
      \int_{-\infty}^{+\infty}h(r)\,\left[\frac{\Ga'}{\Ga}(\tfrac{1}{2}
      +ir)+\frac{\Ga'}{\Ga}(1+ir)\right]\ dr \\
     &\phantom{{=}+2^{\om(m)}}+2\sum_{n=1}^\infty\frac{\La(n)}{n}\,
      \hat h(2\log n)-\sum_{\substack{p|m\\p\text{ prime}}}\sum_{k=0
      }^\infty\frac{\log p}{p^k}\,\hat h(2k\log p)\biggr\}\ .\\
\end{split}
\end{equation*}
\end{prop}

Our main concern will now be to set up the trace formula for $\Ga_0
(m)$ that includes a Hecke operator $T_p$, where $p\ndiv m$. On the 
spectral side, a proof of the trace formula requires to know the 
eigenvalues of $T_p$ not only on a basis for $\cC_m$, but also on 
Eisenstein series. We therefore first establish

\begin{lemma}
\label{Eiseneig}
Let $E_{1/v}(z,s)$, $v|m$, be the non-holomorphic Eisenstein series
$(\ref{Eisendef})$ associated with the representatives $\frac{1}{v}\in 
F(m)$ of inequivalent cusps for $\Ga_0(m)$. Then these are eigenfunctions
of $T_p$, $p\ndiv m$, with eigenvalues $p^{s-\frac{1}{2}}+p^{\frac{1}
{2}-s}$, i.e.,
\begin{equation*}
T_p E_{1/v}(z,s)=\left(p^{s-\frac{1}{2}}+p^{\frac{1}{2}-s}\right)\,
E_{1/v}(z,s)\ .
\end{equation*}
\end{lemma}

\begin{proof}
We first recall \cite[Thm.6.3.3]{Venkov}, which states that
\begin{equation*}
T_p E_{1/v}(z,s)=\sum_{u|m}H_{vu}(s,p)\,E_{1/u}(z,s)\ ,
\end{equation*}  
where the unknowns $H_{vu}(s,p)$ can be obtained as coefficients of 
$y^s$ in a Fourier expansion of $T_p E_{1/v}(\si_u z,s)$. According 
to the decomposition (\ref{decomp}) one obtains  for $\re s>1$
\begin{equation}
\label{TpEisen}
\sqrt{p}\,T_p E_{1/v}(\si_u z,s)=\sum_{\ga\in\Ga_v\backslash M_p(m)}
\left[\im\left(\si_v^{-1}\ga\si_u z\right)\right]^s =\sum_{\tau\in
\Ga_\infty\backslash\si_v^{-1}M_p(m)\si_u}\left[\im\left(\tau z\right)
\right]^s\ .
\end{equation}
As in the case of the Fourier expansion of $E_{1/v}(\si_u z,s)$, see
\cite[sec.3.4]{Iwaniec}, the coefficient of $y^s$ in (\ref{TpEisen}) 
derives from the elements $\tau_\infty=(\begin{smallmatrix}*&*\\0&*
\end{smallmatrix})\in\si_v^{-1}M_p(m)\si_u$. Now consider  $\ga =
\si_v\tau_\infty\si_u^{-1}\in M_p(m)$. Since $\tau_\infty$ fixes
$\infty$, one obtains that $\ga\frac{1}{u}=\frac{1}{v}$. In Lemma 
\ref{ineqcusps}, we found that the $\Ga_0(m)$-inequivalent cusps $\frac{1}{v}$, 
$v|m$, are also inequivalent with respect to $M_p(m)$. Hence $\si_v^{-1}
M_p(m)\si_u\neq\emptyset$ iff $u=v$. This implies that $H_{vu}
(s,p)$ are the entries of a diagonal matrix, and hence the Eisenstein 
series are eigenfunctions of $T_p$. 

We now determine the eigenvalues $H_{vv}(s,p)$. A direct computation 
shows that $\tau_\infty\in\si_v^{-1}M_p(m)\si_v$ is of the form
\begin{equation*}
\tau_\infty=\begin{pmatrix}a+bv&b\frac{v}{m}\\0&d-bv\end{pmatrix}
\ \ \ \ \mbox{with}\ \ \ \ \begin{pmatrix}a&b\\mc&d\end{pmatrix}
\in M_p(m)\ .
\end{equation*}
Thus, one either has (i) $a+bv=p$ and $d-bv=1$, or (ii) $a+bv=1$ 
and $d-bv=p$. These conditions, together with the relation $ad-mbc=p$,
imply that $b(\pm v(p-1)-bv^2 -mc)=0$. Apart from $b=0$, one therefore 
gets the equations $bv+c\frac{m}{v}=\pm (p-1)$. Since $v$ and $\frac{m}
{v}$ are coprime, the solutions $b,c$ are such that $b$ is in a unique 
class mod $\frac{m}{v}$. Furthermore,
\begin{equation*}
\begin{pmatrix}1&N\\0&1\end{pmatrix}\begin{pmatrix}p&b\frac{v}{m}\\
0&1\end{pmatrix}=\begin{pmatrix}p&b\frac{v}{m}+N\\0&1\end{pmatrix}
\ \ \ \ \mbox{and}\ \ \ \ \begin{pmatrix}1&N\\0&1\end{pmatrix}
\begin{pmatrix}1&b\frac{v}{m}\\0&p\end{pmatrix}=\begin{pmatrix}1&b
\frac{v}{m}+Np\\0&p\end{pmatrix}\ .
\end{equation*}
Therefore, there is one equivalence class mod $\Ga_\infty$ of
elements $\tau_\infty\in\si_v^{-1}M_p(m)\si_v$ in the case (i),
whereas there are $p$ classes in the case (ii). Moreover, since
\begin{equation*}
\left[\im\left(\tau_\infty z\right)\right]^s =\begin{cases}
p^s\,y^s&\text{in case (i)},\\p^{-s}\,y^s&\text{in case (ii)},
\end{cases}
\end{equation*}
we conclude that the coefficient of $y^s$ on the right-hand side
of (\ref{TpEisen}) is given by $p^s +p^{1-s}$. Division by $\sqrt{p}$
then yields the eigenvalue of $E_{1/v}(z,s)$.
\end{proof}

\begin{theorem}
\label{STFTpGa0m}
Let $\mu_k =r_k^2+\frac{1}{4}$ run through all eigenvalues of the
hyperbolic Laplacian on $L^2(X_m)$, counted with multiplicities.
Fix a prime $p\ndiv m$ and denote by $t_k(p)$ the eigenvalues 
of $T_p$ on a Hecke basis $\{g_k:\ k\in\nz_0\}$ of $\kz\oplus\cC_m$. 
Then
\begin{equation*}
\begin{split}
\sum_{k=0}^\infty t_k(p)\,h(r_k)= 
     &\sum_{\substack{\{\ga\}_{\Ga_0(m)} \\
     \ga\,\text{ellip.}}}\frac{1}{m_\ga 
      \sqrt{4p-t^2}}\int_{-\infty}^{+\infty}h(r)\,\frac{e^{
      -2r\arcsin\sqrt{1-\frac{t^2}{4p}}}}{1+e^{-2\pi r}}\ dr \\
     &+\frac{1}{\sqrt{p}}\sum_{\substack{\{\ga\}_{\Ga_0(m)},\,
      t\neq p+1 \\\ga\,\text{hyperb.}}}\frac{\arcosh(|\log\ve_\ga|)\,\hat
     h\left(2\arcosh( 
      \tfrac{t}{2\sqrt{p}})\right)}{2\sinh\left(\arcosh(\tfrac{t}{2
      \sqrt{p}})\right)}\\
     &+2^{\om(m)}\biggl\{2\,\hat h(\log p)\biggl[\log\pi+\log(p-1)
      -\frac{\log X(p-1)}{p-1}\biggr] \\
     &\phantom{{=}+2^{\om(m)}}-\tfrac{1}{2} h(0)+\int_{\log p}^\infty
      \hat h(u)\,\frac{e^{\frac{u}{2}}+e^{-\frac{u}{2}}}{e^{\frac{u}{2}}
      -e^{-\frac{u}{2}}+p^{\frac{1}{2}}-p^{-\frac{1}{2}}}\ du \\
     &\phantom{{=}+2^{\om(m)}}-\frac{1}{2\pi}\int_{-\infty}^{+\infty}
      h(r)\,\left[p^{ir}+p^{-ir}\right]\,\frac{\Ga'}{\Ga}(\tfrac{1}{2}
      +ir)\ dr \\
     &\phantom{{=}+2^{\om(m)}}+2\sum_{n=1}^\infty\frac{\La(n)}{n}\,
      \left[\hat h(2\log n -\log p)+\hat h(2\log n+\log p)\right]\\
     &\phantom{{=}+2^{\om(m)}}-\sum_{\substack{q|m \\ q\,\text{prime}}}\sum_{k=0}^{
      \infty}\frac{\log q}{q^k}\,\left[\hat h(2k\log q -\log p)+\hat 
      h(2k\log q+\log p)\right]\bigg\}\ .
\end{split}
\end{equation*}
Here $\displaystyle{X(n)=\prod_{k\!\!\mod n}(k,n)}$. In the elliptic case, $m_\ga$
denotes the order of the centraliser $Z_{\Ga_0(m)}(\ga)=\qz(\ga)\cap
\Ga_0(m)$, and in the hyperbolic case $\ve_\ga$ is a generator of 
$Z_{\Ga_0(m)}(\ga)$.
\end{theorem}

\begin{proof}
To start with, we recall some well-known facts about the trace formula 
for the non-cocompact Fuchsian group $\Ga_0(m)$, compare \cite{Hejhal2,
Iwaniec,Venkov}. One first defines the integral operator
\begin{equation*}
Lg(z)=\int_{\cF_m}K(z,w)\,g(w)\ d\mu(w)\ ,\ \ \ \ 
g\in\cC_m\ ,
\end{equation*}
whose kernel is constructed from a point-pair invariant $k(z,w)=\phi(
\frac{|z-w|^2}{\im z\im w})$, where $\phi$ is an arbitrary function in 
$C_0^2(\rz)$, and
\begin{equation*}
K(z,w)=\sum_{\ga\in\Ga_0(m)/\{\pm E_2\}}k(\ga z,w)\ .
\end{equation*}
Since $\phi$ has compact support, $K(\cdot,w)\in L^2(X_m)$ for any
fixed $w\in\cH$. According to the spectral resolution of the Laplacian, 
one then obtains the spectral expansion
\begin{equation}
\label{Kspecexp}
K(z,w)=\sum_{k=0}^\infty h(r_k)\,g_k(z)\,\overline{g_k}(w)+\frac{1}
{4\pi}\sum_{v|m}\int_{-\infty}^{+\infty}h(r)\,E_{1/v}(z,\tfrac{1}{2}+ir)
\,E_{1/v}(w,\tfrac{1}{2}-ir)\ dr\ ,
\end{equation}
see e.g.\ \cite[Thm.7.4]{Iwaniec}. Here $\{g_k:\ k\in\nz_0\}$ is a 
Hecke basis for $\kz\oplus\cC_m$, and $E_{1/v}(z,s)$ are the Eisenstein 
series (\ref{Eisendef}) associated with $\frac{1}{v}\in F(m)$. The 
function $h(r)$ is related to the point-pair invariant via the Selberg 
transforms \cite[(3.1)]
{Selberg1}
\begin{align}
\label{Seltrans}
Q(t)       &=\int_t^\infty\frac{\phi(u)}{\sqrt{u-t}}\ du \ ,\nonumber\\
\phi(u)    &=-\frac{1}{\pi}\int_u^\infty\frac{Q'(t)}{\sqrt{t-u}}\ dt\ ,\\
\hat h (u) &=Q\left(e^u+e^{-u}-2\right)\ .\nonumber
\end{align}
One furthermore defines
\begin{align*}
\label{Hdef}
H(z,w)  &:=\frac{1}{4\pi}\sum_{v|m}\int_{-\infty}^{+\infty}h(r)\,
         E_{1/v}(z,\tfrac{1}{2}+ir)\,E_{1/v}(w,\tfrac{1}{2}-ir)\ dr \\
K_0(z,w)&:=K(z,w)-H(z,w)\ ,
\end{align*} 
and then proceeds to calculate $\tr K_0$. This finally yields
Proposition \ref{STFGam}. Instead of this, we are here going to consider 
the operator $T_p L$, $p\ndiv m$, with kernel
\begin{equation}
\label{Kpgeom}
K_p(z,w)=\frac{1}{\sqrt{p}}\sum_{j=0}^{p-1}K\left(\frac{z+j}{p},w
\right)+\frac{1}{\sqrt{p}}\,K\left(pz,w\right)=\frac{1}{\sqrt{p}}
\sum_{\ga\in M_p(m)/\{\pm E_2\}}k(\ga z,w)\ ,
\end{equation}
where the last equality follows from (\ref{decomp}). The spectral
expansion of $K_p(z,w)$ can be derived from (\ref{Kspecexp}) and 
Lemma \ref{Eiseneig}, together with the choice of $\{g_k:\ k\in
\nz_0\}$ as a Hecke basis,
\begin{equation}
\label{Kpspec}
\begin{split}
K_p(z,w)=&\sum_{k=0}^\infty h(r_k)\,t_k(p)\,g_k(z)\,\overline{g_k}
          (w) \\
         &+\frac{1}{4\pi}\sum_{v|m}\int_{-\infty}^{+\infty}h(r)\,
          \left[p^{ir}+p^{-ir}\right]\,E_{1/v}(z,\tfrac{1}{2}+ir)
          \,E_{1/v}(w,\tfrac{1}{2}-ir)\ dr\ .
\end{split}
\end{equation}
In order to compute $\tr K_0$, we follow the usual procedure and
first express $K_p (z,z)-T_p H(z,z)$ in terms of the spectral
expansion (\ref{Kpspec}). An integration over the fundamental domain
$\cF_m$ yields the left-hand side of the trace formula in Theorem
\ref{STFTpGa0m}. Then the expansion (\ref{Kpgeom}) will be used
to integrate $K_p (z,z)-T_p H(z,z)$ over the truncated fundamental
domain $\cF_m^Y$, and finally the limit $Y\rto\infty$ will be taken.
Here $\cF_m^Y$ is the fundamental domain $\cF_m$ with
cuspidal regions $\si_v P^Y$ removed for $v|m$, where 
$P^Y :=\{z=x+iy\in\cH:\ -\frac{1}{2}\leq x\leq\frac{1}{2},\ y\geq Y\}$.
We thus have to compute the right-hand side of
\begin{equation}
\label{pretrace}
\begin{split}
\sum_{k=0}^\infty t_k(p)\,h(r_k)=
  &\lim_{Y\rto\infty}\Biggl\{\int_{\cF_m^Y}K_p(z,z)\ d\mu(z)-\frac{1}
   {4\pi}\sum_{v|m}\int_{-\infty}^{+\infty}h(r)\,\left[p^{ir}+p^{-ir}
   \right] \\
  &\phantom{{=}\lim_{Y\rto\infty}\Biggl\{\int_{\cF_m^Y}K_p(z,z)\ 
   d\mu(z)-}\int_{\cF_m^Y}
   E_{1/v}(z,\tfrac{1}{2}+ir)\,E_{1/v}(z,\tfrac{1}{2}-ir)\ d\mu(z)\,dr
   \Biggr\}\ .
\end{split}
\end{equation}

As a first step, we compute the contribution of $T_p H(z,z)$.
To this end we recall the Maa\ss -Selberg relation
\begin{equation}
\begin{split}
\label{MaassSel}
\int_{\cF_m^Y}E_{1/v}(z,\tfrac{1}{2}+ir)\,E_{1/v}(z,\tfrac{1}{2}-ir)
\ d\mu(z)= 
  &2\log Y-\sum_{u|m}\Phi'_{vu}\left(\tfrac{1}{2}+ir\right)
   \Phi^{-1}_{uv}\left(\tfrac{1}{2}+ir\right)\\
  &+\frac{Y^{2ir}}{2ir}\,\Phi_{vv}\left(\tfrac{1}{2}-ir\right)-
   \frac{Y^{-2ir}}{2ir}\,\Phi_{vv}\left(\tfrac{1}{2}+ir\right)+o(1)\ ,
\end{split}
\end{equation}
as $Y\rto\infty$, compare \cite[Prop.6.8.]{Iwaniec}. Here $\Phi(s)=
(\Phi_{uv}(s))$ denotes the scattering matrix for $\Ga_0(m)$, whose 
explicit form can be found in \cite[ch.11,Lem.4.6]{Hejhal2}:
\begin{align*}
\Phi_{vu}(s)&=\vp(s)\,\prod_{p|(v,u)(\frac{m}{v},\frac{m}{u})}
             \frac{p-1}{p^{2s}-1}\prod_{p|(u,\frac{m}{v})(v,\frac{m}
             {u})}\frac{p^s-p^{1-s}}{p^{2s -1}}\ , \\
\vp(s)      &=\sqrt{\pi}\,\frac{\Ga(s-\frac{1}{2})}{\Ga(s)}\,\frac{
             \ze(2s-1)}{\ze(s)}\ .            
\end{align*}
In order to perform the integration of $T_p H(z,z)$ over $\cF_m^Y$, 
see (\ref{pretrace}), we now multiply (\ref{MaassSel}) with 
$\frac{h(r)}{4\pi}\,[p^{ir}+p^{-ir}]$, integrate over $r\in\rz$
and sum over the $2^{\om (m)}$ divisors $v|m$. The first term on the 
right-hand side of (\ref{MaassSel}) then contributes
\begin{equation}
\label{logY}
\frac{2^{\om (m)}}{2\pi}\log Y\int_{-\infty}^{+\infty}h(r)\,\left[
p^{ir}+p^{-ir}\right]\ dr=2^{\om (m)+1}\,\hat h(\log p)\,\log Y\ .
\end{equation}
For the second term, we imitate the computation in \cite[pp.537]
{Hejhal2} and use the known result for $m=1$ as found in \cite[(11.10)]
{HejhalDuke}. This yields
\begin{equation}
\label{Eisencont1}
\begin{split}
-\frac{1}{4\pi}\int_{-\infty}^{+\infty}h(r)\,
  &\left[p^{ir}+p^{-ir}\right]\,\tr\Phi'\Phi^{-1}\left(\tfrac{1}{2}
   +ir\right)\ dr =\\
  &-2^{\om(m)}\biggl\{2\,\hat h(\log p)\,\log\pi-\frac{1}{2\pi}
   \int_{-\infty}^{+\infty}h(r)\,\left[p^{ir}+p^{-ir}\right]\,
   \frac{\Ga'}{\Ga}(\tfrac{1}{2}+ir)\ dr \\
  &\phantom{{=}-2^{\om(m)}}+2\sum_{n=1}^\infty\frac{\La(n)}{n}\,
   \left[\hat h(2\log n -\log p)+\hat h(2\log n+\log p)\right]\\
  &\phantom{{=}-2^{\om(m)}}-\sum_{\substack{q|m \\
  q\,\text{prime}}}\sum_{k=0}^{ 
   \infty}\frac{\log q}{q^k}\,\left[\hat h(2k\log q -\log p)+\hat 
   h(2k\log q+\log p)\right]\bigg\}\ .
\end{split}
\end{equation}
In order to treat the contributions from (\ref{MaassSel}) still 
remaining, we follow \cite[p.155]{Iwaniec} and find
\begin{equation}
\label{Eisencont2}
\begin{split}
\frac{1}{8\pi i}\int_{-\infty}^{+\infty}\frac{h(r)}{r}\,\left[p^{ir}
+p^{-ir}\right]\,
  &\left[Y^{2ir}\,\tr\Phi(\tfrac{1}{2}-ir)-Y^{-2ir}\,\tr\Phi(
   \tfrac{1}{2}+ir)\right]\ dr \\
  &=\tfrac{1}{2}\,\tr\Phi(\tfrac{1}{2})\,h(0) +o(1) \\
  &=2^{\om(m)-1}\,h(0)+o(1)\ .
\end{split}
\end{equation}
This completes the spectral side of the calculation.

Next we integrate $K_p(z,z)$, expressed as in (\ref{Kpgeom}), over
$\cF_m^Y$. First notice that $M_p(m)$ contains only
elliptic and hyperbolic elements. One can therefore proceed as in
the case of the cocompact Fuchsian group $\cO^1$, see the right-hand
side of (\ref{Tptrgeom}). Again, the centraliser $Z_{\Ga_0(m)}(\ga)$ 
of an element $\ga\in M_p(m)$ is given by $\qz(\ga)\cap\Ga_0(m)$. In 
case $\qz(\ga)$ is a field, this is the group of units of norm one in 
an order of the quadratic field $\qz(\ga)\cong\qz(\sqrt{D})$. The only 
exception occurs for $\pm\tr(\ga)=p+1$, where $D=1$ and $\qz(\ga)$ is 
not a field. In this case $Z_{\Ga_0(m)}(\ga)$ only consists of the 
identity, since $\pm\tr(\ga)=p+1$ and $\nm(\ga)=p$ implies that $\ga$ 
has two distinct fixed points on $\qz\cup\{\infty\}$. If there existed  
a non-trivial $\si\in Z_{\Ga_0(m)}(\ga)$, its fixed points would have 
to be $\Ga_0(m)$-equivalent to the ones of $\ga$, i.e.\ they would be 
two distinct cusps of $\Ga_0(m)$, which is impossible. In the 
non-exceptional cases, $\pm\tr(\ga)\neq p+1$, one can proceed as after 
(\ref{Tptrgeom}), thereby already performing the limit $Y\rto\infty$. 
This yields the first two terms on the right-hand side of the trace 
formula in Theorem \ref{STFTpGa0m}.

The exceptional cases require a separate treatment. Their contribution
to the trace formula, before taking the limit $Y\rto\infty$, reads
\begin{equation}
\label{IYint} 
I(Y):=\frac{1}{\sqrt{p}}\sum_{\substack{\ga\in M_p(m) \\ \tr(\ga)=p+1}}
\int_{\cF_m^Y} k(\ga z,z)\ d\mu(z)=\frac{1}{\sqrt{p}}\sum_{\ga_{v,n}}
\int_{\widetilde\cH^Y}k(\ga_{v,n}z,z)\ d\mu(z)\ , 
\end{equation}
where the sum over $\ga_{v,n}$ extends over the representatives of
the exceptional conjugacy classes as described in Lemma \ref{exclemma}.
The domain of integration is given by
\begin{equation}
\label{Htrunc}
\begin{split}
{\widetilde\cH}^Y :=\bigcup_{\al\in\Ga_0(m)}\al\cF_m^Y 
   &=\cH\setminus\bigcup_{u|m}\bigcup_{\al\in\Ga_0(m)}\al\si_u P^Y
     \\
   &=\cH\setminus\bigcup_{u|m}\bigcup_{\tau\in\Ga_\infty\backslash
     \si_u^{-1}\Ga_0(m)\si_u}\si_u\tau P^Y_\infty\ ,
\end{split}
\end{equation}
where $P^Y_\infty :=\{z\in\cH:\ \im z\geq Y\}=\Ga_\infty P^Y$. The 
limit $Y\rto\infty$ cannot be taken before evaluating the integral on 
the right-hand side of (\ref{IYint}), since the domain ${\widetilde
\cH}^Y$ approaches $\cH$ in that limit and the integral over $\cH$ 
diverges. One is therefore forced to evaluate the integral with the 
truncation of the domain present.

In order to proceed further we recall the double coset decomposition
\begin{equation*}
\si_u^{-1}\Ga_0(m)\si_u =\Om_\infty\cup\bigcup_{c>0}\bigcup_{d\!\!\!\!\!
\pmod{mc}}\Om_{d/mc}\ ,
\end{equation*}
see \cite[Thm.2.7]{Iwaniec}. Here the set $\Om_\infty$ is such that 
$\Ga_\infty\backslash\Om_\infty$ consists of a single class with 
representative $\om_\infty =(\begin{smallmatrix}1&*\\0&1\end{smallmatrix})
\in\si_u^{-1}\Ga_0(m)\si_u$. Moreover, $c$ and $d$ run over numbers such that 
$\si_u^{-1}\Ga_0(m)\si_u$ contains $\om_{d/mc}=(\begin{smallmatrix}a&b
\\mc&d\end{smallmatrix})$, and $\Ga_\infty\backslash\Om_{d/mc}=\om_{d/mc}
\Ga_\infty$. Thus
\begin{equation}
\label{horocycles}
\bigcup_{\tau\in\Ga_\infty\backslash\si_u^{-1}\Ga_0(m)\si_u}\tau 
P^Y_\infty =P_\infty^Y\cup\bigcup_{c>0}\bigcup_{d\!\!\!\!\!\pmod{mc}}
\om_{d/mc}P_\infty^Y\ .
\end{equation}
We now exploit the fact that $\al\partial P^Y_\infty$, $\al=
(\begin{smallmatrix}A&B\\C&D\end{smallmatrix})\in SL_2(\rz)$, $C\neq 0$, 
is a horocycle of radius $\frac{1}{2C^2 Y}$ touching $\partial\cH$
at $\frac{A}{C}$, and that $\al^{-1}$ maps $\frac{A}{C}$ to $\infty$. 
In our situation, $\om_{d/mc}P^Y_\infty$ hence is a disc of radius 
$\frac{1}{2m^2c^2 Y}$ touching the real axis in the point $\frac{a}
{mc}$, $(a,mc)=1$, which is mapped to $\infty$ by $\om_{d/mc}^{-1}$. 
This implies that, letting $\om_{d/mc}$ run over the set as specified 
in (\ref{horocycles}), the points $\si_u(\frac{a}{mc})$ are comprised 
of all rational cusps which are $\Ga_0(m)$-equivalent to $\frac{1}{u}$.
According to (\ref{Htrunc}), the domain ${\widetilde\cH}^Y$ therefore
consists of the hyperbolic plane with horocyclic neighbourhoods of all 
cusps removed. 

Special attention has to be devoted to the horocycles touching the real 
axis in the fixed points of $\ga_{v,n}$. The latter are given by 
$-\frac{n}{p-1-nv}$ and $\frac{1}{v}$. In the integral on the right-hand
side of (\ref{IYint}), we change variables from $z$ to
$w=\si_v^{-1}z$. Then $\ga_{v,n}$ gets conjugated to the element
$\si_v^{-1}\ga_{v,n} 
\si_v$, which has fixed points $-\frac{v}{m}\frac{n}{p-1}$ and 
$\infty$.  
We thus notice that the horocycle at the fixed point $-\frac{v}{m}
\frac{n}{p-1}$ corresponds to $\om_{d/mc}P^Y_\infty$ with $a=\frac{vn}
{(vn,mp-m)}$ and $mc=\frac{m(p-1)}{(vn,mp-m)}$. To simplify the integral 
further we now shift the variable $w$ by $\frac{v}{m}\frac{n}{p-1}$, so 
that the fixed points of the transformation $\tilde\ga_{v,n}$ defined 
by a conjugation of $\si_v^{-1}\ga_{v,n}\si_v$ with the shift turn out 
to be $0$ and $\infty$. The matrix representing $\tilde\ga_{v,n}$ is 
hence diagonal. From the knowledge of $\tr(\tilde\ga_{v,n})=p+1$ and 
$\nm(\tilde\ga_{v,n})=p$ one furthermore concludes that one can choose 
$\tilde\ga_{v,n}=(\begin{smallmatrix}p&0\\0&1\end{smallmatrix})$. We 
then denote the shifted domain of integration by $\widehat\cH^Y$ and 
introduce the abbreviation $T:=p^{\frac{1}{2}}-p^{-\frac{1}{2}}$, so that
\begin{equation}
\label{IYint1}
I(Y)=\frac{1}{\sqrt{p}}\sum_{\ga_{v,n}}\int_{\widehat\cH^Y}\phi\left(
\frac{|\tilde\ga_{v,n}z-z|^2}{\im\tilde\ga_{v,n}z\,\im z}\right)\,
d\mu(z)=\frac{1}{\sqrt{p}}\sum_{\ga_{v,n}}\int_{\widehat\cH^Y}\phi
\left(T^2\left(1+\frac{x^2}{y^2}\right)\right)\,\frac{dx\,dy}{y^2}\ .
\end{equation}
By construction, the domain $\widehat\cH^Y$ consists of the hyperbolic
plane $\cH$ where discs, $U_{a/mc}^Y:=\om_{d/mc}P^Y_\infty+\frac{v}{m}
\frac{n}{p-1}$ of specified radii and touching the boundary $\partial\cH$ 
at the cusps $\frac{a}{mc}$, are removed. 

It now turns out that any divergence of $I(Y)$ in the limit $Y\rto
\infty$ comes from the contributions of the cusps that are fixed points 
of the transformation appearing under the integral. In (\ref{IYint1}) 
these cusps are $0$ and $\infty$. To see the finiteness in the other 
cusps, we consider the integral over any of the discs $U_\xi^Y$ for $\xi
\in\qz\setminus\{0\}$,
\begin{equation}
\label{intbound}
\left|\int_{U_\xi^Y}\phi\left(T^2\left(1+\frac{x^2}{y^2}\right)\right)\,
\frac{dx\,dy}{y^2}\right|\leq\mbox{meas}\,\left(U_\xi^Y\right)\cdot
\sup_{z\in U_\xi^Y}\left|y^{-2}\phi\left(T^2\left(1+\frac{x^2}{y^2}\right)
\right)\right|\ ,
\end{equation}   
where $\mbox{meas}\,(U_\xi^Y)$ denotes Lebesgue measure. Since $\phi\in 
C_0^2(\rz)$ is supposed to be compactly supported, $\phi (u)=0$ for 
$|u|>K$ with some sufficiently large constant $K>0$. In particular,
we choose $K>T^2$. Then $\phi(T^2(1+\frac{x^2}{y^2}))=0$ if $y<|x|(
\frac{K}{T^2}-1)^{-\frac{1}{2}}$. In other words, we obtain the bound 
$y^{-2}\leq x^{-2}(\frac{K}{T^2}-1)$ in the support of $\phi$. If $Y$ is 
large enough, then $x\neq 0$ for any $x+iy\in U_\xi^Y$. Thus the 
right-hand side of (\ref{intbound}) is finite, and vanishes as $Y\rto
\infty$. Consequently,
\begin{equation*}
I(Y)=\frac{1}{\sqrt{p}}\sum_{\ga_{v,n}}\int_{\cH\setminus (U_0^Y\cup
U_\infty^Y)}\phi\left(T^2\left(1+\frac{x^2}{y^2}\right)\right)\,\frac{
dx\,dy}{y^2}+o(1)\ .
\end{equation*}
We remark that the integrals depend on $\ga_{v,n}$ through the radius
of the disc $U_0^Y$, since this is a shift of the disc at the fixed
point $\frac{v}{m}\frac{n}{p-1}$ of $\ga_{v,n}$. Therefore, the radius
of $U_0^Y$ is 
\begin{equation}
\label{Radius}
R_{v,n}=\frac{1}{2m^2c^2Y}=\frac{(vn,mp-m)^2}{2m^2(p-1)^2Y}\ .  
\end{equation}
We now introduce polar coordinates $r\in\rz_+$ and $\theta\in [0,\pi]$,
$z=x+iy=r\,e^{i\theta}$, and obtain
\begin{align*}
\int_{\cH\setminus (U_0^Y\cup U_\infty^Y)}\phi\left(T^2\left(1+
\frac{x^2}{y^2}\right)\right)\,\frac{dx\,dy}{y^2}
   &=\int_0^\pi\int_{2R_{v,n}\sin\theta}^{\frac{Y}{\sin\theta}}\phi
     \left(\frac{T^2}{\sin^2\theta}\right)\,\frac{1}{r\,\sin^2\theta} 
     \ dr\,d\theta \\
   &=\int_0^\pi\log\left(\frac{Y}{2R_{v,n}\sin^2\theta}\right)\phi
     \left(\frac{T^2}{\sin^2\theta}\right)\,\frac{1}{\sin^2\theta}\ 
     d\theta \\
   &=\frac{1}{T}\int_{T^2}^\infty\log\left(\frac{Yu}{2R_{v,n}T^2}\right)
     \,\frac{\phi(u)}{\sqrt{u-T^2}}\ du\ .
\end{align*}
In the last line the variable $u:=\frac{T^2}{\sin^2\theta}$ has been 
introduced. In order to proceed further, we recall (\ref{Seltrans})
and find 
\begin{equation*}
\int_{T^2}^\infty\frac{\phi(u)}{\sqrt{u-T^2}}\ du=Q\left(p+p^{-1}-2
\right)=\hat h(\log p)\ .
\end{equation*}
Moreover, using (\ref{Seltrans}) and changing the order of the 
integrations one observes
\begin{equation*}
\begin{split}
\int_{T^2}^\infty\frac{\log u}{\sqrt{u-T^2}}\,\phi(u)\ du
  &=-\frac{1}{\pi}\int_{T^2}^\infty\frac{\log u}{\sqrt{u-T^2}}
   \int_u^\infty\frac{Q'(t)}{\sqrt{t-u}}\ dt\,du \\
  &=-\frac{1}{\pi}\int_{T^2}^\infty Q'(t)\int_{T^2}^t\frac{\log u}
    {\sqrt{u-T^2}\sqrt{t-u}}\ du\,dt\ .
\end{split}
\end{equation*}
The integral over $u$ can be calculated to yield $2\pi\,[\log(T+
\sqrt{t})-\log 2]$, so that an integration by parts yields
\begin{equation*}
\begin{split}
\int_{T^2}^\infty\frac{\log u}{\sqrt{u-T^2}}\,\phi(u)\ du
  &=2\log T\,\hat h(\log p)+\int_{T^2}^\infty\frac{Q(t)}{T+\sqrt{t}}\ 
    \frac{dt}{\sqrt{t}} \\
  &=2\log T\,\hat h(\log p)+\int_{\log p}^\infty\hat h(u)\,\frac{e^{
    \frac{u}{2}}+e^{-\frac{u}{2}}}{e^{\frac{u}{2}}-e^{-\frac{u}{2}}+
    p^{\frac{1}{2}}-p^{-\frac{1}{2}}}\ du\ ,
\end{split}
\end{equation*}
where the last line was obtained through the change of variables $t=
e^u+e^{-u}-2$. Collecting the above results one obtains
\begin{equation*}
\begin{split}
\frac{1}{\sqrt{p}}\int_{\cH\setminus (U_0^Y\cup U_\infty^Y)}\phi\left(
T^2\left(1+\frac{x^2}{y^2}\right)\right)\,\frac{dx\,dy}{y^2}=
  &\frac{\hat h(\log p)}{p-1}\,\log\left(\frac{Y}{2R_{v,n}}\right) \\
  &+\frac{1}{p-1}\int_{\log p}^\infty\hat h(u)\,\frac{e^{\frac{u}{2}}+
   e^{-\frac{u}{2}}}{e^{\frac{u}{2}}-e^{-\frac{u}{2}}+p^{\frac{1}{2}}-
   p^{-\frac{1}{2}}}\ du\ .
\end{split}
\end{equation*}
In order to compute $I(Y)$ we still have to perform the sum over the 
$2^{\om(m)}(p-1)$ representatives $\ga_{v,n}$, which is non-trivial
only for the term containing $R_{v,n}$. According to (\ref{Radius}) and 
Lemma \ref{exclemma} this sum reads
\begin{equation*}
\begin{split}
\sum_{\ga_{v,n}}\log\left(\frac{Y}{2R_{v,n}}\right)
  &=2^{\om(m)+1}(p-1)\,\log Y +2\sum_{\ga_{v,n}}\log\left(\frac{m(p-1)}
    {(vn,mp-m)}\right) \\
  &=2^{\om(m)+1}(p-1)\,\log Y +2\sum_{v|m}\sum_{k=0}^{p-2}\log\left(
    \frac{p-1}{(k,p-1)}\right) \\
  &=2^{\om(m)+1}(p-1)\,\log Y +2^{\om(m)+1}\left[(p-1)\,\log(p-1)-
    \log X(p-1)\right]\ ,
\end{split}
\end{equation*}
where the function $X(n)=\prod\limits_{k\!\!\mod n}(k,n)$ has been 
introduced. Collecting all contributions to $I(Y)$, one first of
all observes that the only term that diverges in the limit $Y\rto
\infty$ is identical to (\ref{logY}). 

According to (\ref{pretrace}) the right-hand side of the trace formula 
follows, if we add the contributions of the elliptic and of the 
non-exceptional hyperbolic conjugacy classes, add $I(Y)$, and subtract 
the contributions (\ref{logY}), (\ref{Eisencont1}) and (\ref{Eisencont2}) 
of the Eisenstein series. The divergent parts therefore cancel, and the
limit $Y\rto\infty$ can be performed. One thus obtains the right-hand
side of the trace formula in Theorem \ref{STFTpGa0m}, however, for
a restricted class of test functions $h$. Since up to now we
required the point-pair invariant $\phi$ to be twice differentiable 
and compactly supported, the same properties hold for $\hat h$.
Then $h$ is smooth and $h(r)=O(r^{-2})$ as $|r|\rto\infty$, see 
\cite[ch.I,Prop.4.1]{Hejhal1}. An extension to the class of test 
functions introduced at the beginning of this section is possible 
since all sums and integrals appearing in the trace formula converge 
absolutely with these test functions. The extension goes via an 
approximation argument as in \cite[pp.32]{Hejhal1} and \cite[Thm.13.8]
{Hejhal2}. This is a standard step and will not be reproduced here, see 
also \cite[Thm.10.2]{Iwaniec}.
\end{proof}

\section{Main theorems}
\label{sec5}

We will need to consider trace formulae for several groups, and to 
distinguish them we let
\[\sum_{g_k\in S}h(r_k) \]
denote the sum over a Hecke basis $\{g_k\}$ of $S$, where $\la_k=r_k^2+
\frac14$ is the Laplace eigenvalue of $g_k$ and $h$ is a test function 
as specified at the beginning of Section \ref{sec4}.

\begin{theorem}
\label{mainthm1}
Let $\cO$ be a maximal order in an indefinite rational quaternion 
division algebra with discriminant $d$. Then the positive Laplace 
eigenvalues, including multiplicities, on $X_{\cO}$ coincide with
the Laplace spectrum on Maa\ss -newforms for the Hecke congruence
group $\Ga_0(d)$.
\end{theorem}

\begin{proof}
We already showed in \cite[Lem.5.4]{BJ1} that all eigenvalues of $-\De$
on $L^2(X_\cO)$ also occur as eigenvalues of $-\De$ on $L^2(X_d)$.
It thus suffices to show that
\begin{equation}
\label{traceid1}
\sum_{g_k\in L^2(X_{\cO})}h(r_k)=\sum_{g_k\in\kz\oplus\cC_d^{new}}h(r_k)
\end{equation}
for an arbitrary test function $h$.

In order to prove (\ref{traceid1}), we will use the right-hand sides of 
the trace formulae in Propositions~\ref{STFcomp} and \ref{STFGam}. To be 
able to apply this procedure to the right-hand side of (\ref{traceid1}), 
we observe that (\ref{newformsum}) implies
\begin{equation}
\label{traceid2}
\sum_{g_k\in\cC_d^{new}}h(r_k)=\sum_{m|d}\be\!\left(\tfrac d m\right)
\sum_{g_k\in\cC_m}h(r_k).
\end{equation}
We now apply the first identity of Lemma \ref{arithid} to obtain
$$
\sum_{m|d}\be\!\left(\tfrac d m\right)\sum_{g_k\in\kz\oplus\cC_m}h(r_k)
=h(r_0)+\sum_{m|d}\be\!\left(\tfrac d m\right)\sum_{g_k\in\cC_m}h(r_k)
=\sum_{g_k\in\kz\oplus\cC_d^{new}}h(r_k)\ .
$$
Using Propositions~\ref{STFcomp} and \ref{STFGam}, we therefore have to 
show that
\begin{equation}
\label{traceid3}
\sum_{g_k\in L^2(X_{\cO})}h(r_k)=\sum_{m|d}\be\!\left(\tfrac d m\right)
\sum_{g_k\in\kz\oplus\cC_m}h(r_k)\ .
\end{equation}

For the terms corresponding to the identity on the right-hand sides of
the trace formulae, one has to show that
\[ A_\cO=\sum_{m|d} \be\!\left(\tfrac d m\right) A_m. \]
This, however, was already proved in \cite{BJ1} using (\ref{area1}) and
(\ref{area2}).

In the elliptic and hyperbolic terms, only the numbers $E'(t,1,\Ga)$ 
depend on the group $\Ga$. Hence, that these terms are identical is 
equivalent to
\[ E'(t,1,\cO^1)=\sum_{m|d} \be\!\left(\tfrac d m\right) 
E'\bigl(t,1,\Ga_0(m)\bigr) \]
for all traces $t$. But this is exactly Lemma~\ref{primitive}.

Now it only remains to prove that the parabolic terms vanish; but this 
follows from Lemma~\ref{arithid} since these are of the  form
\[ \sum_{m|d} \be\!\left(\tfrac d m\right) 2^{\om(m)}
\Bigl(C+\sum_{p|m}f(p)\Bigr),\]
where $C$ is a constant (i.e.\ independent of $m$) and $f$ only depends 
on $p$.
\end{proof}

We define $V_\la^{new}$ to be the subspace of $\cC_{d}^{new}$ with
Laplace eigenvalue $\lambda$, and $W_{\lambda}$ to be the corresponding 
subspace of $L^{2}(X_{\cO})$. The statement of Theorem~\ref{mainthm1} is 
exactly $ \dim V_{\lambda}^{new}=\dim W_{\lambda}$ for all $\lambda$. 
But more generally, we will now prove that the traces of the Hecke 
operators $T_p$ and $\widetilde T_p$, restricted to $V_\la^{new}$ and 
$W_{\lambda}$ respectively, coincide when $p\ndiv d$.

\begin{theorem}
\label{mainthm2}
Let $\cO$ be a maximal order in an indefinite rational quaternion 
division algebra with discriminant $d$. If $p\ndiv d$, then the traces
of $T_{p}$ on $V_{\lambda}^{new}$ and $\widetilde{T}_{p}$ on
$W_{\lambda}$ coincide for all $\lambda$.
\end{theorem}

\begin{proof}
Let $\la=r^2+\frac14>0$ be a Laplace eigenvalue that occurs in the 
newform spectrum of $\Ga_0(d)$. Choose a Hecke basis for $\cC_d^{new}$ 
such that $\{g_1, \ldots, g_N\}$ span $V_\la^{new}$. According to 
Theorem~\ref{mainthm1}, then one can choose a Hecke basis $\{\vp_1,
\ldots,\vp_N\}$ of $W_\la$. Let $t_k(p)$ be the $p$-th Hecke eigenvalue 
of $g_k$ and $\tilde{t}_k(p)$ the $p$-th Hecke eigenvalue of $\vp_k$. We
have to show that 
\begin{equation}
\label{heckeid2}
\sum_{k=1}^N\tilde t_k(p)=\sum_{k=1}^N t_k(p)\ 
\end{equation}
for any prime $p$ with $p\ndiv d$.
In order to prove this identity, we will make use of the trace 
formulae for Hecke operators contained in Proposition~\ref{STFTpO} and 
Theorem~\ref{STFTpGa0m}, respectively, because it suffices to prove that
\begin{equation}
\label{hecketraceid}
\sum_{g_k\in L^2(X_{\cO})}\tilde{t}_k(p)\,h(r_k)=\sum_{g_k\in\kz\oplus
\cC_d^{new}}t_k(p)\,h(r_k)
\end{equation}
for arbitrary test functions $h$. Since $p\ndiv d$, the Hecke operators 
$T_p$ on $\cC_m$ are defined the same way for all $m|d$. We can therefore
employ the same strategy as in the proof of Theorem~\ref{mainthm1}, see 
(\ref{traceid3}), according to which (\ref{hecketraceid}) is equivalent
to
$$
\sum_{g_k\in L^2(X_\cO)}\tilde t_k(p)\,h(r_k)=\sum_{m|d}\be\!\left(
\tfrac{d}{m}\right)\sum_{g_k\in\kz\oplus\cC_m}t_k(p)\,h(r_k)\ .
$$

There are no terms with trace $t=p+1$ occurring in the trace formula for
the left-hand side, since $\cO$ is an order in a division algebra.
The fact that the terms with $t=p+1$ as well as the contributions from
the Eisenstein series to the trace formula for the right-hand side vanish 
follows from Lemma~\ref{arithid}, since these are of the form 
\[ \sum_{m|d} \be\!\left(\tfrac d m\right) 2^{\om(m)}
\Bigl(C+\sum_{\substack{q|m\\q\,\text{prime}}}f(q)\Bigr)\ ,\]
where $C$ is a constant (i.e.\ independent of $m$) and $f$ only depends 
on $q$. 

When comparing the elliptic and the remaining hyperbolic terms, 
one has to be a little more careful than in the proof of 
Theorem~\ref{mainthm1}. First we divide the terms for a fixed trace into 
separate sums corresponding to the different orders $B$ optimally embedded
into $\cO$ and $\Gamma_{0}(m)$, respectively. Then we observe that the 
centraliser, and hence $m_{\ga}$ and $\ve_{\ga}$, only depends on $B$, 
since it is precisely the elements in $B$ with norm equal to $1$. Hence 
the identity of the elliptic and the remaining hyperbolic terms follows 
from (\ref{sumid1}), and thus (\ref{heckeid2}) is established.
\end{proof}

\begin{cor}
\label{cor1}
If $\la>0$ is a Laplace eigenvalue on $L^2 (X_\cO)$ with $\dim W_\la =1$,
then $\Th(W_\la)=V_\la^{new}$.
\end{cor}

\begin{proof}
Let $\vp$ be the element of a Hecke basis for $L^2 (X_\cO)$ that spans
$W_\la$. Then $\Th(\vp)\in\cC_d$ with Laplace eigenvalue $\la$, see
\cite[Prop.5.1]{BJ1}. Moreover, according to \cite[Prop.6.1]{BJ1}
$\widetilde T_p\vp=\tilde t(p)\vp$ implies that $T_p\Th(\vp)=\tilde t(p) 
\Th(\vp)$. If $g$ is the element of a Hecke basis for $\cC_d$ that
spans $V_\la^{new}$, Theorem~\ref{mainthm2} says that $\Th(\vp)$ and 
$g$ have the same Hecke eigenvalues for all $p\ndiv d$. The corollary
then follows from the non-holomorphic analogue of \cite[Lem.20]{AL}.
\end{proof}
 
As a final remark let us mention that for maximal orders $\cO$ the 
Laplace spectrum of $L^2(X_\cO)$ is conjectured to be simple, see
e.g.\ \cite{BJ1}. According to Corollary \ref{cor1}, this would imply 
that $\Th(L_0^2(X_\cO))=\cC_d^{new}$. Following Theorem \ref{mainthm1},
the simplicity of the Laplace spectrum of $L^2(X_\cO)$ is equivalent
to a simple Laplace spectrum on $\cC_d^{new}$. The latter has long been
conjectured, see e.g.\ \cite{Cartier}, and found confirmation in 
extensive numerical calculations of eigenvalues for some groups, see 
e.g.\ \cite{Cartier,Hejhal4,Steil}.


\begin{thebibliography}{10}

\bibitem{AL}
{\sc A.~O.~L.~Atkin and J.~Lehner}, {\em Hecke operators on 
{$\Gamma_0(m)$}}, Math. Ann., 185 (1970), pp.~134--160.

\bibitem{BJ1}
{\sc J.~Bolte and S.~Johansson}, {\em Theta-lifts of Maa\ss\ waveforms},
to appear in Emerging applications of number theory, D.~A.~Hejhal, 
F.~Chung, J.~Friedman, M.~C.~Gutzwiller and A.~Odlyzko, eds., 
Springer-Verlag, New York.

\bibitem{Cartier}
{\sc P.~Cartier}, {\em Some Numerical Computations Relating to Automorphic
Functions}, in Computers in Number Theory, A.~O.~L.~Atkin and J.~B.~Birch, 
eds., Academic Press, London, 1971.

\bibitem{HejhalDuke}
{\sc D.~A.~Hejhal}, {\em The Selberg trace formula and the Riemann zeta 
function}, Duke Math. J., 43 (1976), pp.~441-482.

\bibitem{Hejhal1}
\leavevmode\vrule height 2pt depth -1.6pt width 23pt, {\em The Selberg 
Trace Formula for $PSL(2,\mathbb{R})$ vol. 1}, Lecture Notes in 
Mathematics 548, Springer-Verlag, Berlin-Heidelberg-New York, 1976.

\bibitem{Hejhal2}
\leavevmode\vrule height 2pt depth -1.6pt width 23pt, {\em The Selberg 
Trace Formula for $PSL(2,\mathbb{R})$ vol. 2}, Lecture Notes in 
Mathematics 1001, Springer-Verlag, Berlin-Heidelberg-New York, 1983.

\bibitem{Hejhal3}
\leavevmode\vrule height 2pt depth -1.6pt width 23pt, {\em A classical 
approach to a well-known spectral correspondence on quaternion groups}, 
in Number Theory, New York 1983-84, D.~Chudnovsky, G.~Chudnovsky, 
H.~Cohn, and M.~Nathanson, eds., Lecture Notes in Mathematics 1135,
Springer-Verlag, Berlin-Heidelberg-New York, 1985.

\bibitem{Hejhal4}
\leavevmode\vrule height 2pt depth -1.6pt width 23pt, {\em Eigenvalues
of the Laplacian for Hecke Triangle Groups}, Memoirs of the Amer.\ 
Math.\ Soc., vol.\ 97, No.\ 469, Amer.\ Math.\ Soc., Providence, Rhode
Island, 1992.

\bibitem{Iwaniec}
{\sc H.~Iwaniec}, {\em Introduction to the Spectral Theory of Automorphic
Forms}, Biblioteca de la Revista Matem\'{a}tica Iberoamericana, Madrid,
1995.

\bibitem{JL}
{\sc H.~Jacquet and R.~Langlands}, {\em Automorphic Forms on {GL}(2)}, 
Lecture Notes in Mathematics 114, Springer-Verlag, Berlin-Heidelberg-New 
York, 1970.

\bibitem{Stefan}
{\sc S.~Johansson}, {\em Genera of arithmetic {F}uchsian groups},
to appear in Acta Arith.

\bibitem{Miyake}
{\sc T.~Miyake}, {\em Modular Forms},
Springer-Verlag, Berlin-Heidelberg-New York, 1989.

\bibitem{Schneider}
{\sc V.~Schneider}, {\em Die elliptischen {F}ixpunkte zu {M}odulgruppen 
in {Q}uaternionenschief\-k\"orpern}, Math. Ann., 217 (1975), pp.~29--45.

\bibitem{Selberg1}
{\sc A.~Selberg}, {\em Harmonic analysis and discontinous groups in 
weakly symmetric Riemannian spaces with applications to Dirichlet 
series}, J. Indian Math. Soc., 20 (1956), pp.~47--87.

\bibitem{Steil}
{\sc G.~Steil}, {\em Eigenvalues of the Laplacian and of the Hecke
operators for $PSL(2,\gz)$}, DESY-report 94-028, Hamburg, 1994. 

\bibitem{Venkov}
{\sc A.~B.~Venkov}, {\em Spectral Theory of Automorphic Functions}, Proc.
Steklov Math. Inst., 153 (1981).

\bibitem{Vigneras}
{\sc M.-F.~Vigneras}, {\em Arithm\'etique des Alg\`ebres de Quaternions},
Lecture Notes in Math. 800, Springer-Verlag, Berlin-Heidelberg-New York,
1980.

\end{thebibliography}
\end{document}